\documentclass[a4paper]{article}
\usepackage[utf8]{inputenc}
\usepackage[english]{babel}
\usepackage{authblk}
\usepackage{graphicx}
\usepackage{scrtime}
\usepackage{caption}
\usepackage[colorlinks]{hyperref}
\usepackage{etoolbox}

\usepackage{amsmath,amssymb}
\usepackage{algorithm,algpseudocode}
\usepackage[color=yellow]{todonotes}

\RequirePackage[authoryear]{natbib}

\newtheorem{theorem}{Theorem}[section]

\newtheorem{definition}{Definition}[section]

\newtheorem{lemma}[theorem]{Lemma}

\newtheorem{remark}[theorem]{Remark}

\usepackage{parskip}
\def\keywordname{Keywords}%
\def\keywords#1{\ifx#1\empty\else\def\@keywords{\par\addvspace{10pt}{\keywordfont{\bfseries\keywordname:} #1\par}}\fi}%
\def\@keywords{}%

\begin{document}

\title{Casimir preserving stochastic Lie-Poisson integrators
}


\author[1]{Erwin Luesink \thanks{e.luesink@utwente.nl}}
\author[1]{Sagy Ephrati \thanks{s.r.ephrati@utwente.nl}}
\author[1,2]{Paolo Cifani \thanks{p.cifani@utwente.nl}}
\author[1,3]{Bernard Geurts \thanks{b.j.geurts@utwente.nl}}

\affil[1]{Multiscale Modelling and Simulation, Department of Applied Mathematics, Faculty EEMCS, University of Twente, PO Box 217, 7500 AE Enschede, The Netherlands}
\affil[2]{Gran Sasso Science Institute, Viale F. Crispi 7, 67100 L'Aquila, Italy}
\affil[3]{Multiscale Physics, Center for Computational Energy Research, Department of Applied Physics, Eindhoven University of Technology, PO Box 513, 5600 MB Eindhoven, The Netherlands}

\date{}

\maketitle



\begin{abstract}
Casimir preserving integrators for stochastic Lie-Poisson equations with Stratonovich noise are developed extending Runge-Kutta Munthe-Kaas methods.  The underlying Lie-Poisson structure is preserved along stochastic trajectories. A related stochastic differential equation on the Lie algebra is derived. The solution of this differential equation updates the evolution of the Lie-Poisson dynamics by means of the exponential map. The constructed numerical method conserves Casimir-invariants exactly, which is important for long time integration. This is illustrated numerically for the case of the stochastic heavy top and the stochastic sine-Euler equations.
\end{abstract}

\section{Introduction}
\label{intro}
Many problems in physics and chemistry are described by Hamiltonian models. The most familiar Hamiltonian model is the canonical one, described in terms of momenta and positions. The phase space associated to this kind of model is even dimensional and the dynamics preserves the Hamiltonian itself as well as the canonical symplectic form. A famous example of a canonical Hamiltonian system is the planetary $n$-body problem, see \cite{wisdom1991symplectic}. The equations for the $n$-body problem cannot be solved analytically in general and numerical integration is necessary. However, to avoid the collapse or divergence of the planetary orbits, one requires numerical methods that respect the symplectic structure of the $n$-body problem. 

It often happens that canonical Hamiltonian systems are formulated on the cotangent bundle of a Lie group $G$ and have a $G$-invariant Hamiltonian. Examples of such systems include the rigid body, the heavy top, certain special discretisations of two dimensional ideal hydrodynamics and many problems in quantum mechanics. The $G$-invariance leads to a differentiable symmetry, which by Noether's theorem implies an associated conservation law. The symplectic structure in this setting is replaced by a Lie-Poisson structure. This Lie-Poisson structure is degenerate on certain functions known as Casimirs. In the context of Lie-Poisson (LP) equations, it is the goal of structure-preserving integration to preserve these Casimirs.

In this paper we develop geometric time integration that explicitly incorporates this underlying mathematical structure and preserves important invariants to machine accuracy. The extension to stochastic dynamics that also preserves the underlying geometric structure is the main new result. 

We focus on stochastic systems that have a LP Hamiltonian formulation. Before discussing the stochastic setting, we first review the literature on deterministic LP integration, then we discuss Lie group integration and how LP integration and Lie group integration are related. The review here is far from complete, the monograph \cite{hairer2006geometric} and the review papers \cite{iserles2000lie, marsden2001discrete, celledoni2014introduction} provide a much more complete picture of geometric integration. After the brief review of deterministic LP equations and their numerical integration, we discuss stochastic mechanics, where we focus first on the canonical formulation. We then move on to stochastic LP equations and their numerical integration, which is a topic of recent interest and the main topic of the present work.

The notions of LP structure and LP equations arise from the Marsden-Weinstein reduction theorem, a fundamental result proved in \cite{marsden1974reduction}. Most LP equations cannot be solved analytically and numerical methods are required. A seminal work is \cite{zhong1988lie}, where discrete methods for LP equations were introduced by means of Hamilton-Jacobi theory. This method requires a coordinatisation of the group, which can be difficult to implement. This difficulty was resolved by \cite{channell1991integrators} by reformulating the LP algorithm of \cite{zhong1988lie} in terms of Lie algebra variables. Variational integrators became an alternative to the Hamilton-Jacobi formulation when, in \cite{holm1998euler}, the Lagrangian perspective of Marsden-Weinstein reduction was introduced. The Lagrangian reduction is known as Euler-Poincar\'e reduction. Discrete analogues of Euler-Poincar\'e and LP reduction theory were introduced in \cite{marsden1999discrete} for systems on finite dimensional Lie groups with a $G$-invariant Lagrangian. A review paper that summarises discrete mechanics and variational integrators at that time is \cite{marsden2001discrete}. 

In the pioneering work of \cite{crouch1993numerical}, numerical methods for ordinary differential equations on manifolds were developed. Finite dimensional mechanical systems form an important subclass of differential equations on manifolds. The Crouch-Grossman method, which uses the notion of frames, leads to a class of algorithms with an maximum order of convergence of three. Beyond order three, the analysis of the algorithms becomes very complex. To combat this issue, a different approach was introduced by \cite{munthe1999high}, where numerical methods for differential equations on manifolds were introduced based on Lie groups. These methods are now known as the Runge-Kutta--Munthe-Kaas (RKMK) methods. The review paper \cite{iserles2000lie} summarises most of what was then known about Lie group methods. 

In \cite{engo2001numerical}, the Lie group methods of \cite{munthe1999high} were applied specifically to LP equations. Here the community working on numerical methods for differential equations on manifolds and the community working on discrete mechanics with symmetries intersect. Within this intersection, Lie group methods that preserve the structures associated with mechanics were developed by \cite{bou2007hamilton, bou2009hamilton}. In \cite{bogfjellmo2016high}, symplectic integrators of arbitrarily high order where developed building on the methods of \cite{crouch1993numerical} and \cite{munthe1999high}. A review on Lie group methods can be found in \cite{celledoni2014introduction}. We now move out of the deterministic setting to discuss stochastic mechanics.

Since the work of \cite{bismut1982mecanique}, stochastic Hamiltonian systems have entered as important modelling tools for the analysis of continuous and discrete mechanical systems with uncertainty. In \cite{bismut1982mecanique}, stochastic Hamiltonian systems driven by Brownian motion on symplectic manifolds were introduced. In \cite{lazaro2008stochastic}, the work of Bismut was generalised to manifolds and they showed that the stochastic Hamiltonian systems extremise a stochastic action defined on the space of manifold-valued semimartingales. However, \cite{lazaro2008stochastic} provided a counterexample to the converse statement: ``an extremum of stochastic action satisfies stochastic Hamilton's equations". In \cite{bou2009stochastic}, the focus was restricted to stochastic Hamiltonian systems that are driven by Wiener processes. In this context, \cite{bou2009stochastic} were able to prove almost surely that a curve satisfies stochastic canonical Hamiltonian equations if and only if it extremises a stochastic action. For the numerical integration of stochastic canonical Hamiltonian systems, \cite{bou2009stochastic} introduce stochastic variational integrators in the special case when the Hamiltonian is independent of the momentum variable. In \cite{deng2014high} high order symplectic methods for stochastic canonical Hamiltonian systems were developed. A unifying framework of stochastic discrete variational integrators is developed in \cite{holm2018stochastic}, where the method of \cite{bou2009stochastic} is extended to general Hamiltonian systems. Drift preserving methods for stochastic Hamiltonian systems were introduced in \cite{chen2020drift} and variational integrators for stochastic diffusive Hamiltonian systems were developed in \cite{kraus2021variational} by means of stochastic Lagrange-D'Alembert variational principles. 

In \cite{holm2015variational}, a stochastic variational principle was introduced with the purpose of deriving stochastic Euler-Poincar\'e equations. Stochastic Euler-Poincar\'e equations are equivalent to stochastic LP equations when the Legendre transform is a diffeomorphism, so the stochastic variational principle provides a systematic means to derive stochastic LP equations. The type of stochasticity introduced in \cite{holm2015variational} is called ``stochastic advection by Lie transport (SALT)" and its purpose is to model principally unknown effects influencing transport in fluid problems. In \cite{arnaudon2018noise}, this stochastic variational principle is used to derive equations for finite dimensional mechanical systems. The SALT noise belongs to a class of ``Hamiltonian noise", which does not affect the Poisson bracket. This means that the Casimirs for a system that is perturbed with Hamiltonian noise are the same as for the unperturbed system. This provides the desire for structure-preserving numerical methods for stochastic LP equations.

The structure-preserving numerical integration of stochastic LP equations has recently received attention. In \cite{brehier2023splitting} splitting methods for Stratonovich stochastic LP equations are introduced. These methods are explicit and preserve Casimirs as well as the Poisson map property. We will take a different approach to the numerical solution of Stratonovich stochastic LP equations. Our approach is based on the RKMK methods of \cite{munthe1999high} and \cite{engo2001numerical}, which means that it preserves Casrimir invariants associated with the LP structure. The method has a beneficial property compared to other methods with regards the strong order of convergence for the kind of noise considered in \cite{holm2015variational}.

The work presented here is a first step in creating high-order invariant-preserving RKMK methods for solving stochastic LP systems. The use of the RKMK framework is emphasized here as it opens up the possibility for an extension toward high-order methods in the future, such as the 1.5 and 2.0 strong-Taylor schemes illustrated in \cite{kloeden1992stochastic, holm2018stochastic}. The important step of achieving lower-order invariant-preserving geometric time integration for stochastic LP systems is worked out in detail here. The invariant-preserving property is built into the numerical method and illustrated explicitly for the well-known heavy top problem as well as for the fluid-mechanical sine-Euler model developed in \cite{zeitlin1991finite} which connects the current work to applications in geophysics. Different numerical tests are presented to underpin the practical usefulness of the new time integration. We restrict the illustrations to low-dimensional dynamics - this is not a principal restriction of the new geometric time integration method. In fact, extension to high-performance simulation of fully resolved spatial models of Navier-Stokes type was achieved in \cite{cifani2022casimir} recently.

This paper is structured as follows. In Section \ref{sec:slp} we introduce LP brackets to which the geometric structure of the problem is associated. In Section \ref{sec:slpdynamics} we introduce stochastic LP dynamics by introducing semimartingale Hamiltonians. In Section \ref{sec:ni} we introduce the new numerical stochastic LP integrator. This integrator is able to preserve the Casimirs and upon removing the noise, recovers the results of \cite{engo2001numerical}. In Section \ref{sec:ht} we apply the trapezoidal rule combined with the Munthe-Kaas method (TMK) to the stochastic heavy top and show that the Casimirs are conserved, whereas the trapezoidal rule applied directly to the LP equations fails to exactly conserve the Casimirs. In Section \ref{sec:sEuler}, we provide further illustration of the TMK method for the sine-Euler equations \cite{zeitlin1991finite}. The latter system is characterised by higher-degree (polynomial) invariants, the conservation of which is crucial for proper capturing of the intricate dynamics associated with this system. We will show that the developed stochastic geometric integrator is, in fact, able to preserve such structure. Conclusions are collected in Section \ref{sec:conclusion}.

\section{Lie-Poisson brackets}\label{sec:slp}
We start by recalling some notation and definitions of \cite{marsden2013introduction}. We will denote Lie groups by $G$ and the associated Lie algebra by $\mathfrak{g}$. An LP bracket is a linear Poisson bracket on a vector space. Let $f\in C^\infty(V)$, with $V$ a vector space, and let the duality pairing $\langle\,\cdot\,,\cdot\,\rangle_{V^*\times V}:V^*\times V\to\mathbb{R}$ define the dual $V^*$. One can determine the variational derivative $\delta f/\delta u\in V^*$ as the unique element defined by
\begin{equation}
\delta f(u) = \lim_{\epsilon\to 0} \frac{1}{\epsilon}\big(f(u+\epsilon \delta u)-f(u)\big) = \left\langle \delta u,\frac{\delta f}{\delta u}\right\rangle_{V^*\times V},
\end{equation}
with $u\in V$ and $\delta u\in V$ arbitrary. Vector spaces with linear Poisson brackets are duals of Lie algebras. This follows from the fact that the dual $\mathfrak{g}^*$ of any Lie algebra $(\mathfrak{g},[\,\cdot\,,\,\,\cdot\,])$ carries a linear Poisson bracket
\begin{equation}
\{f,g\}_{\mp}(\mu) := \mp\left\langle\mu,\left[\frac{\delta f}{\delta\mu},\frac{\delta g}{\delta \mu}\right]\right\rangle,
\label{eq:kksbracket}
\end{equation}
where $f,g\in C^\infty(\mathfrak{g}^*)$, $\mu\in\mathfrak{g}^*$. The variational derivatives $\delta f/\delta\mu, \delta g/\delta \mu$ are elements of the dual of the dual Lie algebra $\mathfrak{g}^{**}\simeq\mathfrak{g}$ and $\langle\,\cdot\,,\cdot\,\rangle:\mathfrak{g}^*\times\mathfrak{g}\to\mathbb{R}$ is the nondegenerate pairing between the Lie algebra and its dual. This formulation is valid for both finite and infinite dimensional Lie algebras. LP brackets are degenerate. This means that there exist functions $c\in C^\infty(\mathfrak{g}^*)$ called Casimirs such that $\{f,c\}_{\mp}=0$ for all $f\in C^\infty(\mathfrak{g}^*)$.

Casimirs are elements of the kernel of the LP bracket, which means that they are conserved by LP dynamics. In the adjoint and coadjoint representation theory of the Lie algebra, the Lie bracket can be expressed as 
\begin{equation}
\left\langle \mu,\left[\frac{\delta f}{\delta \mu}, \frac{\delta g}{\delta \mu}\right]\right\rangle = \left\langle \mu, {\rm ad}_{\delta f/\delta \mu}\frac{\delta g}{\delta \mu}\right\rangle = \left\langle {\rm ad}^*_{\delta f/\delta \mu} \mu, \frac{\delta g}{\delta \mu}\right\rangle,
\end{equation} 
where ${\rm ad}:\mathfrak{g}\times\mathfrak{g}\to\mathfrak{g}$ is the adjoint representation of the action of the Lie algebra on itself and ${\rm ad}^*:\mathfrak{g}\times\mathfrak{g}^*\to\mathfrak{g}^*$ is the coadjoint representation of the action of the Lie algebra on its dual. There are also representations of the action of the Lie group on the Lie algebra and its dual, given by ${\rm Ad}:G\times\mathfrak{g}\to\mathfrak{g}$ and ${\rm Ad}^*:G\times\mathfrak{g}^*\to\mathfrak{g}^*$. For matrix Lie groups, ${\rm Ad}_g v = gvg^{-1}$ and $\langle{\rm Ad}^*_{g^{-1}}\mu,v\rangle = \langle \mu, {\rm Ad}_g v\rangle$ defines the coadjoint action.

The set $\mathcal{O}_{\mu_0}:=\{{\rm Ad}^*_{g^{-1}}\mu_0\,|\,g\in G\}$ defined by the LP bracket is called the coadjoint orbit of $\mu_0\in\mathfrak{g}^*$. On coadjoint orbits, the Casimirs are constant. This follows from the fact that the Casimirs are in the kernel of ${\rm ad}^*$ and ${\rm ad}^*$ is the infinitesimal generator for ${\rm Ad}^*$. The coadjoint orbits will be key in the construction of the Casimir-preserving numerical integrator. For deterministic systems, there is a large class of different LP integrators that are able to reflect this property numerically, see for instance \cite{zhong1988lie}, \cite{channell1991integrators}, \cite{mclachlan1993explicit}, \cite{reich1994momentum}, \cite{mclachlan1995equivariant}, \cite{engo2001numerical} for different developments of numerical LP integrators. More recently the intersection of isospectral methods and LP methods led to novel LP integrators as shown in \cite{bloch2006isospectral} and \cite{modin2020lie}.

Henceforth, we will only consider finite dimensional Lie algebras. For finite dimensional Lie algebras with dimension $N$, let $e_i$ denote the basis so that a Lie algebra element $\sigma\in\mathfrak{g}$ can be expressed as $\sigma = \sum_{i=1}^N\hat{\sigma}^i e_i$. Let $\varepsilon^i$ denote the induced dual basis by the Frobenius pairing. Then an element $\mu\in\mathfrak{g}^*$ can be expressed as $\mu=\sum_{i=1}^N \hat{\mu}_i\varepsilon^i$. The LP bracket \eqref{eq:kksbracket} can be expressed in terms of the structure constants $C_{ij}^k$ of the Lie algebra $\mathfrak{g}$ as follows
\begin{equation}
\{f,g\}_{\mp}(\mu) := \mp \sum_{i,j,k=1}^N C_{ij}^k \mu_k\frac{\partial f}{\partial \mu_i}\frac{\partial g}{\partial \mu_j},
\label{eq:strucconstbracket}
\end{equation}

The relation between \eqref{eq:kksbracket} and \eqref{eq:strucconstbracket} is expressed in the following lemma.
\begin{lemma}\label{lemma:lpbracket}
For any $\sigma\in\mathfrak{g}$ and any $\mu\in\mathfrak{g}^*$
\begin{equation}
{\rm ad}_\sigma^* \mu = - \sum_{i,j,k=1}^N C_{ij}^k \hat{\mu}_k \hat{\sigma}^j\varepsilon^i =: - J(\mu)\sigma.
\label{eq:ad*Jy}
\end{equation}
\end{lemma} 

Lemma \ref{lemma:lpbracket} establishes a useful relation between the coadjoint representation of the Lie algebra on its dual, the structure constants of the Lie algebra and the skew-symmetric matrix $J$. In the canonical case, the matrix $J$ is the familiar symplectic matrix. The coadjoint representation is important since it can be used to solve LP equations exactly at a formal level. In addition lemma \ref{lemma:lpbracket} defines the linear operator $J$, which generalises the symplectic matrix encountered in canonical Hamiltonian dynamics.

LP brackets arise naturally after reducing the dimension in mechanical systems with symmetry, as is described for the rigid body in \cite{smale1970topologya}, \cite{smale1970topologyb} and in general in \cite{marsden1974reduction}, \cite{marsden2013introduction}, \cite{holm2008geometric2}, \cite{holm2009geometric}. Such mechanical systems can be formulated on Lie groups. When the Hamiltonian is invariant under the left or right action of that Lie group, one can perform symmetry reduction and obtain left or right Lie-Poisson equations as is established in \cite{marsden1974reduction}. Since chirality plays a role, the LP brackets \eqref{eq:kksbracket} and \eqref{eq:strucconstbracket} feature $\mp$, with the minus sign corresponding to the left invariant situation and the plus sign corresponding to the right invariant situation. For $\mu\in\mathfrak{g}^*$, LP equations can be formulated in several equivalent forms. For the sake of compact notation, we use Einstein's summation convention, i.e., the summation will be understood over lower and upper pairs of indices, where the summation runs to $N$, the dimension of the Lie algebra $\mathfrak{g}$.

The LP equations associated to a Hamiltonian $\hslash:\mathfrak{g}^*\to\mathbb{R}$ are given for a general $f\in C^\infty(\mathfrak{g}^*)$ by
\begin{equation}
\begin{aligned}
\frac{d}{dt}f(\mu) &= \mp\left\langle\mu,\left[\frac{\delta f}{\delta \mu},\frac{\delta\hslash}{\delta\mu}\right]\right\rangle\\
&= \mp C_{ij}^k \mu_k \frac{\partial f}{\partial \mu_i}\frac{\partial \hslash}{\partial \mu_j},
\end{aligned}
\end{equation}
which implies that the momentum $\mu$ satisfies the equation
\begin{equation}
\begin{aligned}
\frac{d}{dt}\mu &= \pm {\rm ad}^*_{\delta\hslash/\delta\mu} \mu,\\
&= \mp J(\mu)\frac{\delta\hslash}{\delta \mu},
\end{aligned}
\label{eq:leftlp}
\end{equation}
where ${\rm ad}^*:\mathfrak{g}\times\mathfrak{g}^*\to\mathfrak{g}^*$ is the dual of the adjoint representation ${\rm ad} = [\,\cdot\,,\,\cdot\,]:\mathfrak{g}\times\mathfrak{g}\to\mathfrak{g}$ and we have used Lemma \ref{lemma:lpbracket}.

The left- and right-invariant cases indicate that chirality introduces a sign difference at this stage. However, chirality expresses itself more subtly in several subsequent derivations. To emphasise the (small) differences, we have occasionally split the text into two columns, one column corresponding to the left-invariant case and the other column corresponding to the right-invariant case. At this stage, one can formulate a deterministic LP integrator to solve equation \eqref{eq:leftlp}. This is what \cite{engo2001numerical} did and they proceeded to write a deterministic LP integrator that preserves the Hamiltonian. We will first introduce stochasticity into the LP equations in section \ref{sec:slpdynamics} and formulate an integrator in Section \ref{sec:ni}.

\section{Stochastic Lie-Poisson dynamics}\label{sec:slpdynamics}
Following \cite{protter2005stochastic}, we introduce a filtered probability space given by the quadruplet $(\Omega,\mathcal{F},(\mathcal{F}_t)_{t\geq 0},\mathbb{P})$. Here $\Omega$ is a set, $\mathcal{F}$ is the $\sigma$-algebra, $(\mathcal{F}_t)_{t\geq 0}$ is a right-continuous filtration and $\mathbb{P}$ is the probability measure.

With respect to the filtered probability space, we define a family $W_t^1,\hdots,W_t^M$ of independent, identically distributed Brownian motions. In this section, we will assume that all considered stochastic processes are compatible with the continuous semimartingale $\mathbf{S}_t=(S_t^0, S_t^1, \hdots, S_t^M)$, where $M\leq \dim \mathfrak{g}$. Compatibility is understood in the sense of Definition 2.3 in \cite{street2021semi}, which says that all stochastic processes have Radon-Nikodym derivatives with respect to each other. The symbol $\circ$ will indicate that the stochastic integral is to be understood in the Stratonovich sense, instead of composition of functions for which this symbol is also adopted. The Stratonovich integral is preferred for the development of stochastic LP equations, because Stratonovich processes satisfy the ordinary chain rule, meaning that a Hamiltonian plus semimartingale noise is enough to determine the equations of motion. The It\^o integral can also be used, but this requires a connection on the underlying manifold, see \cite{emery2006two}, \cite{huang2023second}.

Noise will be introduced via a semimartingale Hamiltonian $\hslash_s:\mathfrak{g}^*\to\mathbb{R}$ and a continuous semimartingale $\mathbf{S}_t$. It is possible to define stochastic LP equations with general semimartingales. We will discuss the conditions for the existence and uniqueness of strong solutions for stochastic LP equations with general semimartingales. For our examples in Sections \ref{sec:ht} and \ref{sec:sEuler}, we restrict to stochastic Hamiltonian systems driven by Wiener processes with specific diffusion Hamiltonians, as dictated by the SALT framework. In this case, the semimartingale $\mathbf{S}_t$ takes the explicit form
\begin{equation}
\mathbf{S}_t = (t,W_t^1,\hdots,W_t^M),
\label{def:semibrownian}
\end{equation}
and semimartingale Hamiltonian splits as
\begin{equation}
\hslash_s\circ d\mathbf{S}_t := \hslash\,dt + \sum_{i=1}^M \widetilde{\hslash}_i \circ dW_t^i.
\label{eq:semimartingalehamiltoniangen}
\end{equation}

We associate the Hamiltonian $\hslash:\mathfrak{g}^*\to\mathbb{R}$ with the drift $dt$ such that the deterministic LP equations are recovered upon setting the diffusion Hamiltonians $\widetilde{\hslash}_i:\mathfrak{g}^*\to\mathbb{R}$ to zero. The noise is controlled by diffusion Hamiltonians (or noise Hamiltonians) $\widetilde{\hslash}_i$ which are associated with the diffusions $dW_t^i$.

As in the deterministic case, the stochastic LP equations distinguish a left-invariant version and a right-invariant version. The left- and right-invariant stochastic LP equations are given, respectively, by
\begin{equation}
{\sf d}\mu \mp {\rm ad}^*_{\partial\hslash_s/\partial\mu} \mu \circ d\mathbf{S}_t = 0.
\label{seq:leftlpad*}
\end{equation}

Let the initial datum be $\mu(0) = \mu_0\in\mathfrak{g}^*$. When $\mathbf{S}_t$ is a general semimartinagle, the integral form of equation \eqref{seq:leftlpad*} is
\begin{equation}
\mu(t) = \mu_0 \pm \int_0^t {\rm ad}^*_{\partial\hslash_s/\partial\mu} \mu(s) \circ d\mathbf{S}_s,
\label{seq:lpsemimartingale}
\end{equation}
for which Theorems 6 and 7 in chapter 5 of \cite{protter2005stochastic} can be applied after the It\^o correction to prove that strong solutions to \eqref{seq:lpsemimartingale} exist, are semimartingales themselves and are unique provided that the integrand is functional Lipschitz. The precise definition of a functional Lipschitz operator can be found in \cite{protter2005stochastic}, together with the chain of implications
\begin{equation}
\begin{aligned}
\text{Lipschitz function} &\implies\text{Random Lipschitz function}\\
\text{Random Lipschitz function} &\implies\text{Process Lipschitz operator} \\
\text{Process Lipschitz operator} &\implies\text{Functional Lipschitz operator}.
\end{aligned}
\end{equation}
Hence, a sufficient condition for strong wellposedness of \eqref{seq:lpsemimartingale} is that the semimartingale Hamiltonian is a $C^{2,1}$ function, that is, the semimartingale Hamiltonian is twice differentiable with Lipschitz continuous derivatives. One derivative is required for computing the gradient of the Hamiltonian, a second derivative is required for the It\^o correction and the result needs to be Lipschitz continuous. For weak solutions to \eqref{seq:lpsemimartingale}, one can consider the conditions discussed in \cite{stroock1997multidimensional} and for more general results on conditions for strong wellposedness of solutions we refer to \cite{jacod2006calcul}. 

In case the semimartingale is of the form
\begin{equation}
\mathbf{S}_t = (t,W_t^1,\hdots,W_t^M),
\end{equation} 
the conditions for strong wellposedness of solutions to the stochastic LP equation reduce to the conditions for Stratonovich diffusions, which require the drift to be Lipschitz continuous and the diffusions to be $C^{1,1}$, because one needs to do the It\^o correction. This means that the drift Hamiltonian needs to be in $C^{1,1}$ and the diffusion Hamiltonians need to be in $C^{2,1}$, i.e., twice differentiable with Lipschitz continuous derivatives.

A canonical choice for a noise Hamiltonian $\widetilde{\hslash}_i$ is through coupling of noise to the momentum $\mu\in\mathfrak{g}^*$. This corresponds to the concept of ``stochastic advection by Lie transport" (SALT) that was introduced in \cite{holm2015variational}, \cite{de2020implications} and is also adopted in this paper. For SALT, the semimartingale Hamiltonian in \eqref{eq:semimartingalehamiltoniangen} takes the form
\begin{equation}
\hslash_s \circ d\mathbf{S}_t = \hslash\,dt + \sum_{i=1}^M \beta_i\cdot\mu \circ dW_t^i,
\label{eq:semimartingalehamiltonian}
\end{equation}
with $\beta_i\in\mathfrak{g}$ being the Lie algebra-valued noise coefficient so that $\beta_i\cdot\mu\in\mathbb{R}$ for each $i$. The Lie algebra-valued noise coefficient determines the amplitude of the noise in the direction of each basis vector. Since the diffusion Hamiltonians are smooth in this case, the stochastic LP equation has unique strong solutions provided that the drift Hamiltonian is in $C^{1,1}$. Upon inserting the semimartingale Hamiltonian \eqref{eq:semimartingalehamiltonian} into the LP bracket \eqref{eq:strucconstbracket}, the resulting stochastic LP equation will have linear multiplicative noise:
\begin{equation}
{\sf d}\mu \mp {\rm ad}^*_{\partial \hslash/\partial \mu}\mu \,dt \mp \sum_{i=1}^M {\rm ad}^*_{\beta_i}\mu\circ dW_t^i = 0.
\label{eq:lpsde}
\end{equation} 

The adjoint and coadjoint representation theory for Lie algebras and their dual can be used to solve the LP equations \eqref{seq:leftlpad*} formally. Adjoint representation theory features the linear operators ${\rm Ad}:G\times \mathfrak{g}\to\mathfrak{g}$ and ${\rm ad}:\mathfrak{g}\times\mathfrak{g}\to\mathfrak{g}$ and coadjoint representation theory features the linear operators ${\rm Ad}^*:G\times\mathfrak{g}^*\to\mathfrak{g}^*$ and ${\rm ad}^*:\mathfrak{g}\times\mathfrak{g}^*\to\mathfrak{g}^*$. The operator ${\rm Ad}^*$ is the dual of ${\rm Ad}$ and ${\rm ad}^*$ is the dual of ${\rm ad}$ with respect to the pairing $\langle\,\cdot\,,\,\cdot\,\rangle:\mathfrak{g}^*\times\mathfrak{g}\to\mathbb{R}$. These operators play a fundamental role in the solution of equations on Lie algebras. The following lemma shows the differential relations between the operators, where we have split the text into columns to emphasise the role of chirality.\bigskip

\begin{lemma}
Let $g(t):\mathbb{R}_+ \to G$,  $\sigma(t):\mathbb{R}_+\to\mathfrak{g}$ and $\mu(t):\mathbb{R}_+\to\mathfrak{g}^*$ all depend almost surely continuously on $t$ and be compatible with the semimartingale $\mathbf{S}_t$. Then the following formulas hold
\bigskip

\noindent
\begin{minipage}{.46\textwidth}
\centering
Left-invariant case
\rule{\textwidth}{1pt}
\flushleft
Denote by ${\sf d}\zeta = g^{-1}({\sf d}g)$ the left-invariant vector field,
\begin{equation}
\begin{aligned}
{\sf d}{\rm Ad}_g \sigma &= {\rm Ad}_{g}\left({\sf d}{\sigma} + {\rm ad}_{{\sf d}\zeta}\sigma\right),\\
{\sf d}{\rm Ad}_{g^{-1}}\sigma &= {\rm Ad}_{g^{-1}} {\sf d}{\sigma} - {\rm ad}_{{\sf d}\zeta}{\rm Ad}_{g^{-1}}\sigma,\\
{\sf d}{\rm Ad}^*_{g^{-1}}\mu &= {\rm Ad}^*_{g^{-1}}{\sf d}{\mu} + {\rm ad}^*_{{\sf d}\zeta}{\rm Ad}^*_{g^{-1}}\mu,\\
{\sf d}{\rm Ad}^*_{g}\mu &= {\rm Ad}^*_{g}({\sf d}{\mu} - {\rm ad}^*_{{\sf d}\zeta}\mu).
\end{aligned}
\label{seq:leftcoadorbitrelations}
\end{equation}
\end{minipage}
\begin{minipage}{.05\textwidth}
\hfill
\end{minipage}
\begin{minipage}{.46\textwidth}
\centering
Right-invariant case
\rule{\textwidth}{1pt}
\flushleft
Denote by ${\sf d}\zeta = ({\sf d}g)g^{-1}$ the right-invariant vector field,
\begin{equation}
\begin{aligned}
{\sf d}{\rm Ad}_g \sigma &= {\rm Ad}_{g}{\sf d}{\sigma} + {\rm ad}_{{\sf d}\zeta} {\rm Ad}_g \sigma,\\
{\sf d}{\rm Ad}_{g^{-1}}\sigma &= {\rm Ad}_{g^{-1}}({\sf d}{\sigma} - {\rm ad}_{{\sf d}\zeta}\sigma),\\
{\sf d}{\rm Ad}^*_{g^{-1}}\mu &= {\rm Ad}^*_{g^{-1}}({\sf d}{\mu} + {\rm ad}^*_{{\sf d}\zeta}\mu),\\
{\sf d}{\rm Ad}^*_{g}\mu &= {\rm Ad}^*_{g}{\sf d}{\mu} - {\rm ad}^*_{{\sf d}\zeta}{\rm Ad}^*_{g}\mu.
\end{aligned}
\label{seq:rightcoadorbitrelations}
\end{equation}
\end{minipage} 
\end{lemma}
\medskip

The proof is a direct computation using the relations between the vector field ${\sf d}\zeta$ and the group element $g$ as well as the definitions of the representations ${\rm Ad}, {\rm Ad}^*, {\rm ad}$ and ${\rm ad}^*$. From the fourth equation in \eqref{seq:leftcoadorbitrelations} it follows that the solution to the left-invariant stochastic LP equation is given by
\begin{equation}
\mu(t) = {\rm Ad}^*_{g(t)} \mu_0,
\end{equation}
and from the third equation in \eqref{seq:rightcoadorbitrelations} it follows that the solution to the right-invariant stochastic LP equation is given by
\begin{equation}
\mu(t) = {\rm Ad}^*_{g(t)^{-1}} \mu_0.
\end{equation}
In both cases, the curve $g(t):\mathbb{R}\to G$ is a continuous curve that is the solution to a stochastic differential equation related to the LP equations. The curve $\mu(t):\mathbb{R}_+\to\mathfrak{g}^*$ that solves the left- or right-invariant stochastic LP equation lives on the set $\mathcal{O}_{\mu_0}$ defined by
\begin{equation}
\mathcal{O}_{\mu_0}:=\{{\rm Ad}^*_{g^{-1}}\mu_0\,|\, g\in G\}.
\label{eq:coadjointorbitdef}
\end{equation} 
The set $\mathcal{O}_{\mu_0}$ is the coadjoint orbit generated by the initial condition. Since Casimirs are in the kernel of the operator ${\rm ad}^*$, this implies that Casimirs are constant on coadjoint orbits. To construct a numerical method that preserves the Casimirs exactly, we need to numerically compute the coadjoint orbits exactly. This can be realised by computing a numerical approximation to the group element $g$ that itself is an element of the group $G$. This has been pioneered by \cite{engo2001numerical} for deterministic LP systems. We use that as starting point for the extension toward stochastic problems. 
\bigskip

\section{Numerical integration}\label{sec:ni}
In \cite{engo2001numerical}, implicit numerical methods that preserve Casimirs as well as the Hamiltonian were introduced, based on the Runge-Kutta Munthe-Kaas (RKMK) integrators of \cite{munthe1999high} and the discrete gradient method of \cite{gonzalez1996time} for integration of ordinary differential equations on manifolds. The discrete gradient method is designed for the conservation of first integrals, which can be used to obtain conservation of the Hamiltonian. We will employ similar methods for the integration of stochastic LP equations, but we first introduce the RKMK method in the deterministic setting.

The RKMK method is based on canonical coordinates of the first kind. Following the notation of \cite{bou2009hamilton}, one introduces a map $\tau:\mathfrak{g}\to G$ such that $\tau$ is a local diffeomorphism of neighbourhood of $0\in \mathfrak{g}$ to a neighbourhood of the identity $e\in G$ with $\tau(0)=e$. The map $\tau$ is assumed to be analytic in the neighbourhood of the identity and is also required to satisfy $\tau(X)\tau(-X)=e$ for all elements $X\in\mathfrak{g}$. This means that $\tau$ induces a local chart on $G$ for which left translation can be used to form an atlas. The exponential map $\exp:\mathfrak{g}\to G$ is an important example of a canonical coordinate inducing map. One can construct RKMK methods based on any $\tau:\mathfrak{g}\to G$ by computing the derivative $d\tau$ of $\tau$ and its inverse $d\tau^{-1}$. For a general $\tau$, the RKMK method is defined below.
\begin{definition}
Given RK coefficients $b_i, a_{ij}\in\mathbb{R}$ $(i,j=1,\hdots,s)$, set $c_i=\sum_{j=1}^s a_{ij}$. An $s$-stage RKMK approximant to the differential equation $\dot{g}(t) = g(t)f(t,g(t))$ with initial datum $g(0)=g_0\in G$ is given by 
\begin{equation}
\begin{aligned}
G_k^i &= g_k\tau(h\Theta^i_k),\\
\Theta^i_k &= h \sum_{j=1}^s a_{ij} d\tau^{-1}_{-h\Theta^j_k}f(t_k + c_j h, G^j_k),\quad i=1,\hdots,s,\\
g_{k+1} &= g_k \tau\left(h\sum_{j=1}^s b_j d\tau^{-1}_{-h\Theta_k^j} f(t_k+c_j h, G_k^j)\right).
\end{aligned}
\end{equation}
If $a_{ij}=0$ for $i<j$, then the RKMK method is explicit and it is implicit otherwise. 
\end{definition}
In most cases $d\tau^{-1}$ is expressed by a series expansion that will have to be truncated to be able to use the RKMK method. Theorem 4.7 in \cite{bou2009hamilton} explains that given an $r$th order approximation to the exact exponential map $\tau$ and a RK method of order $p$ with $r\geq p$, then the RKMK method is of order $p$ if the truncation of $d\tau^{-1}$ satisfies $q\geq p-2$.

By Ado's theorem (see \cite{rossmann2006lie}, page 51), we have that every finite dimensional Lie algebra is isomorphic to a matrix Lie algebra. This means that for computational purposes, the exponential map can be expressed as the matrix exponential
\begin{equation}
\exp(X) = \sum_{k=0}^\infty\frac{X^k}{k!}.
\end{equation}
This means that we can compute the differential $d\exp$ and its inverse $d\exp^{-1}$ explicitly, which will be used in the construction for the RKMK method for stochastic LP equations.

By means of the linear representations of Lie algebras and their duals, one can establish the following important relationship.

\begin{lemma}\label{lemma:ad*relation}
For any $\sigma\in\mathfrak{g}$,
\begin{equation}
{\rm Ad}^*_{\exp(\sigma)} = \exp(-{\rm ad}^*_\sigma).
\end{equation}
\end{lemma}
\medskip

\noindent Lemma \ref{lemma:ad*relation} is the coadjoint version of the fundamental relation ${\rm Ad}_{\exp(\sigma)}=\exp({\rm ad}_\sigma)$. The proof is a straightforward calculation, see e.g. \cite{rossmann2006lie}. Lemma \ref{lemma:ad*relation} shows that $-{\rm ad}^*$ is the infinitesimal generator of ${\rm Ad}^*$. It also implies that Casimirs, which are in the kernel of ${\rm ad}^*$ result in the identity operator. Hence the value of the Casimir is constant on coadjoint orbits. We will now construct the RKMK method for stochastic LP equations. We will represent a group element $g=\exp\sigma$ (locally) by a Lie algebra element. As a consequence, the following representation of the solution $\mu(t)$ is obtained. 
\medskip

\noindent
\begin{minipage}{.46\textwidth}
\centering
Left-invariant case
\rule{\textwidth}{1pt}
\flushleft
The solution to the left-invariant LP equation is given by
\begin{equation}
\mu(t) = {\rm Ad}^*_{g}\mu_0.
\label{seq:leftAd*g}
\end{equation}
The next step is to apply the stochastic product rule to \eqref{seq:leftAd*g}. This requires a derivative of the exponential map, since the vector field ${\sf d}\zeta$ is given by ${\sf d}\zeta=g^{-1}{\sf d}{g}$ in this case. Taking the stochastic differential and using the third equation in \eqref{seq:leftcoadorbitrelations} yields
\begin{equation}
\begin{aligned}
{\sf d} \mu &= {\sf d}{\rm Ad}^*_{g} \mu_0\\
&= {\rm ad}_{{\sf d}\zeta}^* {\rm Ad}^*_{g} \mu_0\\
&= {\rm ad}_{{\sf d}\zeta}^* \mu.
\end{aligned}
\end{equation}
\end{minipage}
\begin{minipage}{.05\textwidth}
\hfill
\end{minipage}
\begin{minipage}{.46\textwidth}
\centering
Right-invariant case
\rule{\textwidth}{1pt}
\flushleft
The solution to the right-invariant LP equation is given by
\begin{equation} \tag{23}
\mu(t) = {\rm Ad}^*_{g^{-1}}\mu_0.
\label{seq:rightAd*g}
\end{equation}
The next step is to apply the stochastic product rule to \eqref{seq:rightAd*g}. This requires a derivative of the exponential map, since the vector field ${\sf d}\zeta$ is given by ${\sf d}\zeta={\sf d}{g}g^{-1}$ in this case. Taking the stochastic differential and using the fourth equation in \eqref{seq:rightcoadorbitrelations} yields
\begin{equation} \tag{24}
\begin{aligned}
{\sf d}\mu &= {\sf d}{\rm Ad}^*_{g^{-1}}\mu_0\\
&= -{\rm ad}_{{\sf d}\zeta}^* {\rm Ad}^*_{g^{-1}} \mu_0\\
&= -{\rm ad}_{{\sf d}\zeta}^* \mu
\end{aligned}
\end{equation}
\end{minipage}

\noindent
\begin{minipage}{.46\textwidth}
Since $\mu$ satisfies the stochastic LP equation, we can deduce that
\begin{equation}
{\sf d}\zeta = \frac{\delta \hslash_s}{\delta \mu}\circ d\mathbf{S}_t.
\label{seq:leftxih}
\end{equation}
Now we compute ${\sf d}\zeta$ from its definition ${\sf d}\zeta=g^{-1}{\sf d}{g}$ with $g=\exp(\sigma)$
\begin{equation}
\begin{aligned}
{\sf d}\zeta &= \exp(-\sigma){\sf d}\exp(\sigma)\\
&=  d\exp_{\sigma}({\sf d}{\sigma}).
\end{aligned}
\label{seq:leftxiexp}
\end{equation}
\end{minipage}
\begin{minipage}{.05\textwidth}
\hfill
\end{minipage}
\begin{minipage}{.46\textwidth}
Since $\mu$ satisfies the stochastic LP equation, we can deduce that
\begin{equation} \tag{25}
{\sf d}\zeta = \frac{\delta\hslash_s}{\delta\mu}\circ d\mathbf{S}_t.
\label{seq:rightxih}
\end{equation}
Now we compute $\xi$ from its definition $\xi={\sf d}{g}g^{-1}$ with $g=\exp(\sigma)$
\begin{equation} \tag{26}
\begin{aligned}
{\sf d}\zeta &= {\sf d}\exp(\sigma) \exp(-\sigma)\\
&= d\exp_{-\sigma}({\sf d}{\sigma}).
\end{aligned}
\label{seq:rightxiexp}
\end{equation}
\end{minipage}

\noindent
\begin{minipage}{.46\textwidth}
By using \eqref{seq:leftxih} and \eqref{seq:leftxiexp}, we obtain the stochastic differential equation for $\sigma$, which is given by
\begin{equation} \tag{22}
\begin{aligned}
{\sf d}\sigma &= d\exp_{\sigma}^{-1}\left(\frac{\delta\hslash_s}{\delta\mu}\right)\circ d\mathbf{S}_t.
\end{aligned}
\label{seq:leftsigmaeq}
\end{equation}
\end{minipage}
\begin{minipage}{.05\textwidth}
\hfill
\end{minipage}
\begin{minipage}{.46\textwidth}
By using \eqref{seq:rightxih} and \eqref{seq:rightxiexp}, we obtain the stochastic differential equation for $\sigma$, which is given by
\begin{equation} \tag{27}
\begin{aligned}
{\sf d}\sigma &= d\exp_{-\sigma}^{-1}\left(\frac{\delta\hslash_s}{\delta\mu}\right)\circ d\mathbf{S}_t.
\end{aligned}
\label{seq:rightsigmaeq}
\end{equation}
\end{minipage}
\medskip

The operator in equation \eqref{seq:leftsigmaeq} acting on the variational derivative of the Hamiltonian is defined through its power series expansion as
\begin{equation} \tag{28}
\begin{aligned}
d\exp_{\sigma}^{-1}(v) &= \sum_{k=0}^\infty \frac{B_k^+}{k!}{\rm ad}_\sigma^k(v)\\
&= v + \frac{1}{2}[v,\sigma] + \frac{1}{12}[[v,\sigma],\sigma] + \hdots 
\end{aligned}
\label{eq:bernoullileft}
\end{equation}
where $B_k^+$ are the Bernoulli numbers with the convention $B_1^+= +1/2$ and ${\rm ad}^k_\sigma(v)$ is defined recursively as ${\rm ad}^k_\sigma(v) = [\sigma,{\rm ad}_\sigma^{k-1}(v)]$. Similarly, in the right-invariant case in \eqref{seq:rightsigmaeq}, we have
\begin{equation} \tag{29}
\begin{aligned}
d\exp_{-\sigma}^{-1}(v) &= \sum_{k=0}^\infty \frac{B_k^-}{k!}{\rm ad}_\sigma^k(v)\\
&= v - \frac{1}{2}[v,\sigma] + \frac{1}{12}[[v,\sigma],\sigma] + \hdots
\end{aligned}
\label{eq:bernoulliright}
\end{equation}
where $B_k^-$ are the Bernoulli numbers with the convention $B_1^-=- 1/2$.

Note that the operator $d\exp_{\sigma}^{-1}(\,\cdot\,):\mathfrak{g}\to\mathfrak{g}$ is a linear operator. The iterated commutators are nonlinear in $\sigma$ and linear in the other variable. The SDEs \eqref{seq:leftsigmaeq} and \eqref{seq:rightsigmaeq} can be discretised using any numerical method that is consistent with Stratonovich SDEs, see \cite{kloeden1992stochastic} for instance. As a result, the Casimirs will be preserved by construction. 

This is the essence of the deterministic RKMK method developed in \cite{munthe1999high}. It was shown in \cite{munthe1999high} that one needs as many terms in the expansions \eqref{eq:bernoullileft} and \eqref{eq:bernoulliright} as the order of convergence of the method minus one to preserve the order of convergence of the overall method. Casimirs are guaranteed to be conserved numerically, but to also conserve energy one has to put in some more work. \cite{engo2001numerical} used the discrete derivative method of \cite{gonzalez1996time} together with the RKMK method to obtain energy and Casimir conservation. Energy conserving methods based on the discrete derivative method of \cite{gonzalez1996time} are second order in time, which implies that only the first term in the expansions \eqref{eq:bernoullileft} and \eqref{eq:bernoulliright} is required.

The energy and Casimir conserving method that \cite{engo2001numerical} developed yields a discretisation based on the trapezoidal rule. In terms of the variable $\mu$ on the dual of the Lie algebra, this leads to the natural discretisation for the Stratonovich integral, see e.g., \cite{stratonovich1966new}, \cite{protter2005stochastic}. It is in general not possible to construct arbitrarily high order methods for stochastic differential equations based on solely increments of Wiener processes, as was shown by \cite{clark1980maximum}. One can use iterated integrals to go to higher order, but this is a notoriously difficult problem, see \cite{kloeden1992stochastic}. Thus, in practice only the first term in the expansions \eqref{eq:bernoullileft} and \eqref{eq:bernoulliright} is required. This reduces the operator $d\exp_{\sigma}^{-1}(\,\cdot\,)$ to the identity operator. The SDE for $\sigma$ with initial condition $\sigma(0) = 0$ is then given by
\begin{equation} \tag{30}
\begin{aligned}
{\sf d}\sigma &= \frac{\delta\hslash_s}{\delta\mu}\circ d\mathbf{S}_t\\
&= \frac{\delta \hslash}{\delta\mu}\,dt + \sum_{i=1}^M\frac{\delta\widetilde{\hslash}_i}{\delta\mu}\circ dW_t^i.
\end{aligned}
\label{seq:sigma}
\end{equation}
The SDE \eqref{seq:sigma} can be solved with any appropriate method such that the diffusion term converges towards the Stratonovich integral. Since the Stratonovich integral is the limit of the midpoint approximation to this integral, the Heun method (explicit midpoint) and the implicit midpoint method are natural choices. Note that the SDE \eqref{seq:sigma} depends on the Hamiltonian, which in turn depends on the momentum $\mu$. This is important for the solution procedure that we will now construct.
\medskip

\begin{remark} \label{remark:additivenoise}
Let us consider the case where the noise Hamiltonians are chosen according to SALT. This means that the noise Hamiltonians are linear in the momentum, which implies that the stochastic Lie-Poisson equations have linear multiplicative noise of the Stratonovich type. Equation \eqref{seq:sigma} has additive noise in the SALT case. This facilitates the analysis of stochastic Lie-Poisson equations and implies that the strong order of convergence of the Euler-Maruyama method is one, instead of one half.
\end{remark}
\medskip

In the following we employ the trapezoidal rule to discretise equation \eqref{seq:sigma}. The discretisation of \eqref{seq:sigma} viewed in terms of $\mu$ coincides with the midpoint rule, but when viewed solely in terms of $\sigma$, the discretisation is the trapezoidal rule. The trapezoidal rule coincides with the discrete derivative approach of \cite{gonzalez1996time}, which was designed to conserve integrals of Hamiltonian dynamics. In \cite{engo2001numerical}, the discrete derivative method was used to design a class of energy conserving Lie-Poisson integrators. Here we will use the discrete derivative approach because it corresponds to the midpoint discretisation of the Stratonovich integral and in absence of noise recovers the results of \cite{engo2001numerical}. Applying the discrete derivative of \cite{gonzalez1996time} to both the drift and the diffusion terms in \eqref{seq:sigma} yields the following discretisation
\begin{equation} \tag{31}
\sigma_n = \Delta t \frac{\delta}{\delta\mu}\left(\hslash(\mu_{n+\frac{1}{2}})\right) + \sum_{i=1}^M\Delta W^i_{n} \frac{\delta}{\delta \mu}\left(\widetilde{\hslash}_i(\mu_{n+\frac{1}{2}})\right),
\label{seq:implicit}
\end{equation}
where $\Delta t$ is the time-step size, $\Delta W^i_{n} = W^i_{n+1} - W^i_{n}$ and $\mu_{n+\frac{1}{2}} = \frac{1}{2}(\mu_n + \mu_{n+1})$. Note that by definition the second term in \eqref{seq:implicit} converges to the Stratonovich integral as the time step size tends to zero. The Hamiltonian is updated every time step, which means that the differential equation \eqref{seq:sigma} is solved for only one time step every update. This means that $\sigma_{n-1}$ is always equal to the initial condition, which is zero. Without stochasticity, \eqref{seq:implicit} coincides with the method of \cite{engo2001numerical} that conserves the deterministic Hamiltonian. To find an estimate for $\mu_{n+1}$, we use a type of quasi-Newton method called the chord method to find the root of a function. This function will be equation \eqref{seq:implicit} with $\sigma_n$ subtracted from both sides. We will use Lemma \ref{lemma:ad*relation} to write this function as
\begin{equation} \tag{32}
\begin{aligned}
f(\sigma_n) &= \sigma_n - \Delta t\left(\frac{\delta}{\delta\mu}\hslash\left(\frac{1}{2}\mu_n + \frac{1}{2}\exp(\mp{\rm ad}^*_{\sigma_n})\mu_n\right)\right)\\
&\quad - \sum_{i=1}^M\Delta W^i_n \left(\frac{\delta}{\delta\mu}\widetilde{\hslash}_i\left(\frac{1}{2}\mu_n+ \frac{1}{2}\exp(\mp{\rm ad}^*_{\sigma_n})\mu_n\right)\right).
\end{aligned}
\end{equation}
We now want to determine when $f(\sigma_n)=0$, as this will yield the $\sigma_n$ that is required to go from $\mu_n$ to $\mu_{n+1}$. In the Newton-Raphson method, the Jacobian has to be updated every time step and then inverted
\begin{equation} \tag{33}
\begin{aligned}
\sigma_n^{[k+1]} &= \sigma_n^{[k]} - \left(Df(\sigma_n^{[k]})\right)^{-1} f(\sigma_n^{[k]}).
\end{aligned}
\label{eq:newtonmethod}
\end{equation}
This can be very expensive. Instead we use the chord method, which freezes the Jacobian on the initial condition. Computing the Jacobian and evaluating on $\sigma_n^{[0]}=0$ results in
\begin{equation} \tag{34}
Df(0) = I + \frac{\Delta t}{2}D^2\hslash(\mu_n)J(\mu_n) + \sum_{i=1}^M \frac{\Delta W_n^i}{2} D^2\widetilde{\hslash}_i(\mu_n)J(\mu_n).
\end{equation}
Here $D^2 \hslash$ denotes the Hessian of a functional $\hslash$. Hence, the chord method is given by
\begin{equation} \tag{35}
\begin{aligned}
\sigma_n^{[k+1]} &= \sigma_n^{[k]} - \left(Df(0)\right)^{-1} f(\sigma_n^{[k]}).
\end{aligned}
\label{eq:chordmethod}
\end{equation}
We conclude this section by summarising the stochastic Lie-Poisson integrator based on the trapezoidal rule.
\begin{algorithm}[H]
\caption{Stochastic Lie-Poisson integrator based on TMK}
\begin{algorithmic}
\State $n \gets 0$
\State $\mu \gets \mu_0$
\State $\Delta W \gets \sqrt{\Delta t}\,{\rm randn}(M,N)$
\While{$n \neq N$}
\State $n \gets n+1$
\State $\sigma \gets {\rm Chord\, method}(\mu)$ \quad (see equation \eqref{eq:chordmethod})
\State $g \gets \exp(\sigma)$
\State $\mu \gets {\rm Ad}^*_{g} \mu$
\EndWhile
\end{algorithmic}
\label{alg:IMMK}
\end{algorithm}
\medskip

In section \ref{sec:ht} we will apply the stochastic Lie-Poisson integrator specified in algorithm \ref{alg:IMMK} above to the stochastic heavy top.

\section{Heavy top}\label{sec:ht}
\begin{figure}[h!]
\centering
\includegraphics[width=.6\textwidth]{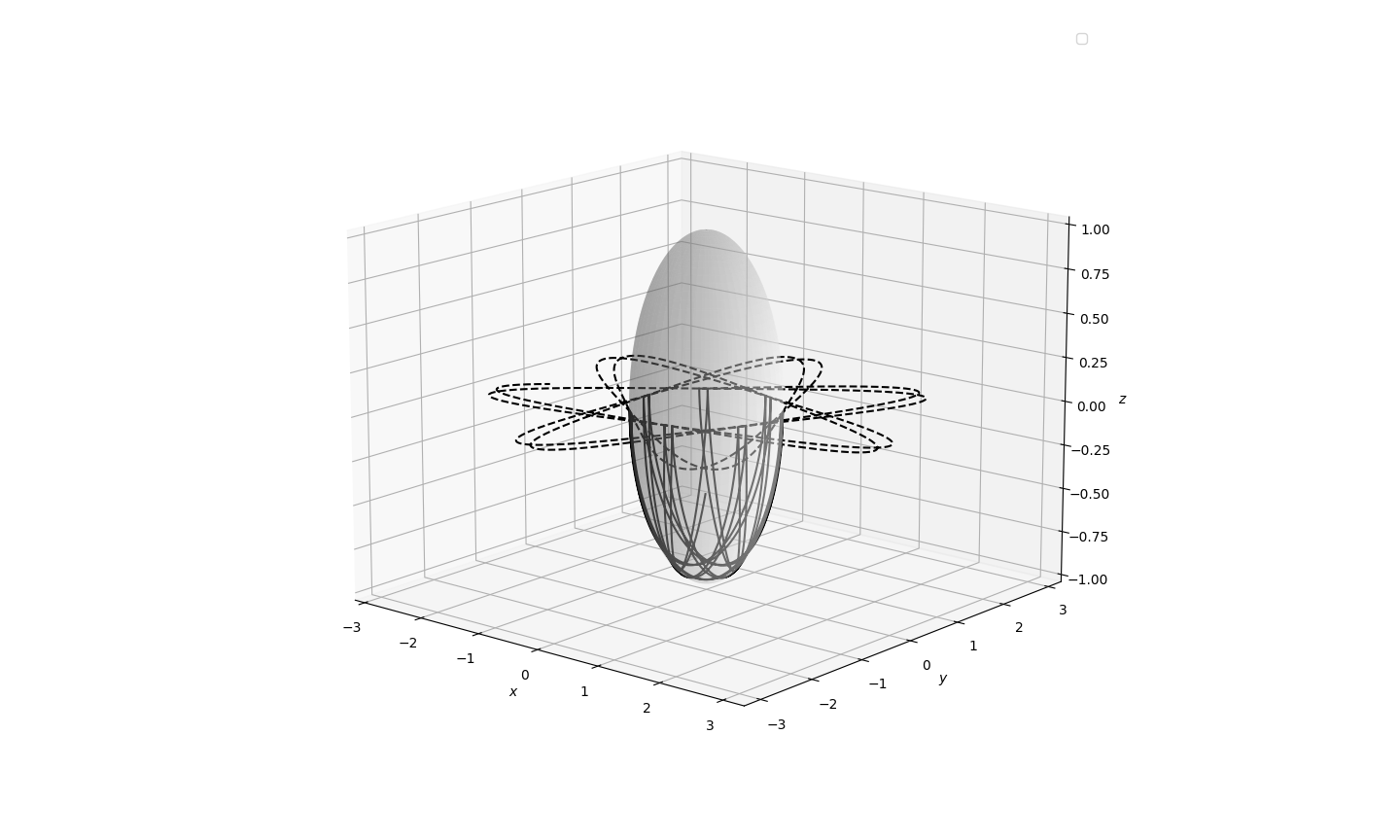}
\caption{A deterministic trajectory of the heavy top with moment of inertia given by $\mathbb{I}={\rm diag}(4,4,1)$ generated with the TMK method. The angular momentum variable (dashed) $\boldsymbol \pi$ is at constant $z$ value in the $xy$-plane. The integrable nature of this particular setting explains the regular periodic orbits.}
\label{fig:detdynamics}
\end{figure}

In this section we will apply the stochastic Lie-Poisson integrator to the stochastic heavy top. The stochastic heavy top is a Lie-Poisson system on the Lie algebra $\mathfrak{se}(3)\subset\mathbb{R}^{4\times 4}$ associated to the special Euclidean group $SE(3)\subset\mathbb{R}^{4\times 4}$. The special Euclidean group is the group of rotations and translations. A representation of $SE(3)$ is as follows
\begin{equation} \tag{36}
SE(3) = \left\{ \left[\begin{matrix}
\mathbf{R} & \mathbf{v}\\
\mathbf{0}^T & 1
\end{matrix}\right]\, \bigg| \, \mathbf{R}\in SO(3)\,\text{ and }\,\mathbf{v}\in\mathbb{R}^3\right\},
\end{equation}
where $\mathbf{R}\in SO(3)\subset\mathbb{R}^{3\times 3}$ is a rotation matrix. The Lie algebra $\mathfrak{se}(3)$ has the associated representation
\begin{equation} \tag{37}
\mathfrak{se}(3) = \left\{ \left[\begin{matrix}
\widehat{\boldsymbol{\xi}} & \mathbf{b}\\
\mathbf{0}^T & 1
\end{matrix}\right]\, \bigg| \, \widehat{\boldsymbol \xi}\in \mathfrak{so}(3)\,\text{ and }\, \mathbf{b}\in\mathbb{R}^3\right\},
\end{equation}
where $\widehat{\boldsymbol \xi}\in\mathfrak{so}(3)\subset\mathbb{R}^{3\times 3}$ is a traceless, skew-symmetric matrix. The matrix exponential applied to an element of $\mathfrak{se}(3)$ yields an element of $SE(3)$. Traceless, skew-symmetric matrices $\widehat{\boldsymbol \xi}\in\mathbb{R}^{3\times 3}$ have three degrees of freedom and can be represented by a vector in $\boldsymbol \xi \in \mathbb{R}^3$ via the hat map isomorphism $\,\widehat{}\,:\mathbb{R}^3\to\mathfrak{so}(3)$. The hat map isomorphism permits us to represent all variables associated with the heavy top in the sequel by vectors in $\mathbb{R}^3$. The semimartingale Hamiltonian for the stochastic heavy top chosen to be
\begin{equation} \tag{38}
\begin{aligned}
\hslash_s(\boldsymbol \pi,\boldsymbol \gamma)\circ d\mathbf{S}_t &= \hslash\,dt + \widetilde{\hslash}\circ dW_t\\
&= (\boldsymbol \pi\cdot \mathbb{I}^{-1}\boldsymbol \pi + \boldsymbol \chi\cdot\boldsymbol \gamma)\, dt + \boldsymbol \alpha\cdot\boldsymbol \pi \circ dW_t.
\end{aligned}
\label{eq:semimarthamht}
\end{equation}
Here $\boldsymbol \pi\in\mathbb{R}^3$ is the angular momentum vector. The moment of inertia tensor is represented by the diagonal matrix $\mathbb{I}={\rm diag}(I_1,I_2,I_3)\in\mathbb{R}^{3\times 3}$. The vector $\boldsymbol \gamma\in\mathbb{R}^3$ is called the gravity vector and tracks the direction of gravity relative to the orientation of the top. The vector $\boldsymbol \chi\in\mathbb{R}^3$ connects the point around which the top rotates towards the centre of gravity of the top. The vector $\boldsymbol \alpha\in\mathbb{R}^3$ is the vector of noise amplitudes. The semimartingale Hamiltonian \eqref{eq:semimarthamht} involves the usual Hamiltonian for the heavy top $\hslash$ and a single noise Hamiltonian $\widetilde{\hslash}$ which couples noise to the momentum. The equations of motions of the stochastic heavy top are the following stochastic Lie-Poisson equations
\begin{equation} \tag{39}
\begin{aligned}
{\sf d}\boldsymbol \pi &= (\boldsymbol \pi \times \boldsymbol \omega + \boldsymbol \chi\times\boldsymbol \gamma)\,dt + \boldsymbol \pi\times \boldsymbol \alpha\circ dW_t,\\
{\sf d}\boldsymbol \gamma &= \boldsymbol \gamma \times \boldsymbol \omega\,dt + \boldsymbol \gamma\times\boldsymbol \alpha\circ dW_t,\\
\boldsymbol \omega &= \mathbb{I}^{-1}\boldsymbol \pi,
\end{aligned}
\label{eq:liepoissonht}
\end{equation} 
with initial conditions $\boldsymbol \pi(0) = \boldsymbol \pi_0\in\mathbb{R}^3$ and $\boldsymbol \gamma(0) = \boldsymbol \gamma_0\in\mathbb{R}^3$. The heavy top is a completely integrable system in a number of situations.  If $I_1=I_2=I_3$, then the heavy top is known as the Euler top, which is completely integrable. If two moments of inertia are the same and center of gravity lies on the symmetry axis, then we deal with a Lagrange top, this is also a completely integrable system as shown in \cite{ratiu1980motion}. Two more integrable cases are known, the heavy top is a completely integrable system if $I_1=I_2=2I_3$ (the Kovalevskaya top) or if $I_1=I_2=4I_3$ (the Goryachev-Chaplygin top). We investigate two cases, in the first case we introduce SALT-type stochasticity to the Goryachev-Chaplygin top and in the second case we focus on a nonintegrable situation. We will compare standard implicit midpoint (IM) integration applied directly to the stochastic LP equations \eqref{eq:liepoissonht} to the integrator introduced above (we will refer to this integrator as Trapezoidal Munthe-Kaas method or TMK method for short). The dynamics of the gravity vector $\boldsymbol \gamma(t)$ is restricted to the sphere with radius $|\boldsymbol \gamma_0|$ as a result of the $\mathfrak{so}(3)$ action, whereas the dynamics of the angular momentum $\boldsymbol \pi(t)$ wanders more freely in $\mathbb{R}^3$. A deterministic trajectory of the symmetric heavy top with moment of inertia $\mathbb{I}={\rm diag}(4,4,1)$ is displayed in figure \ref{fig:detdynamics}.

For stochastic simulations, we first take the Goryachev-Chaplygin top by setting the moment of inertia to be $\mathbb{I}={\rm diag}(4,4,1)$. A second simulation concerns a nonintegrable case with moment of inertia $\mathbb{I}={\rm diag}(4,2,1)$. Realisations of these two stochastic heavy tops are shown in Figures \ref{fig:dynamics1}. The sphere with radius $|\boldsymbol \gamma_0|$ is shaded in grey. For the simulations, we have used time step size $\Delta t = 0.01$ and simulated until time $T = 100$. If one takes step size beyond a certain value $\Delta t \geq \Delta t_c$ (whose precise value depends on which stochastic LP system one is considering, the parameters of said LP system and the realisation of the noise), the properties of the method deteriorate. For simulations of the Goryachev-Chaplygin top, this critical value was empirically found to be typically around $\Delta t_c \approx 0.8$ and for the nonintegrable case $\Delta t_c \approx 0.6$. The time step size $\Delta t=0.01$ is chosen such that in both cases it is an order of magnitude smaller than the critical values.

The initial conditions are $\boldsymbol \pi_0 = \frac{1}{2}\sqrt{2}(-1,1,0)$ and $\boldsymbol \gamma_0 = \frac{1}{2}\sqrt{2}(-1,1,0)$. The centre of gravity lies above the point around which the top rotates, so $\boldsymbol \chi = (0,0,1)$. The noise amplitude is $\boldsymbol \alpha = (0.01,0.02,0.03)$. In Figure \ref{fig:dynamics2}, the moment of inertia is $\mathbb{I}={\rm diag}(4,4,1)$ and in Figure \ref{fig:dynamics1}, the moment of inertia is $\mathbb{I}={\rm diag}(4,2,1)$.

\noindent
\begin{minipage}[t]{.47\textwidth}
\vspace{-2ex}
\includegraphics[width=.99\textwidth]{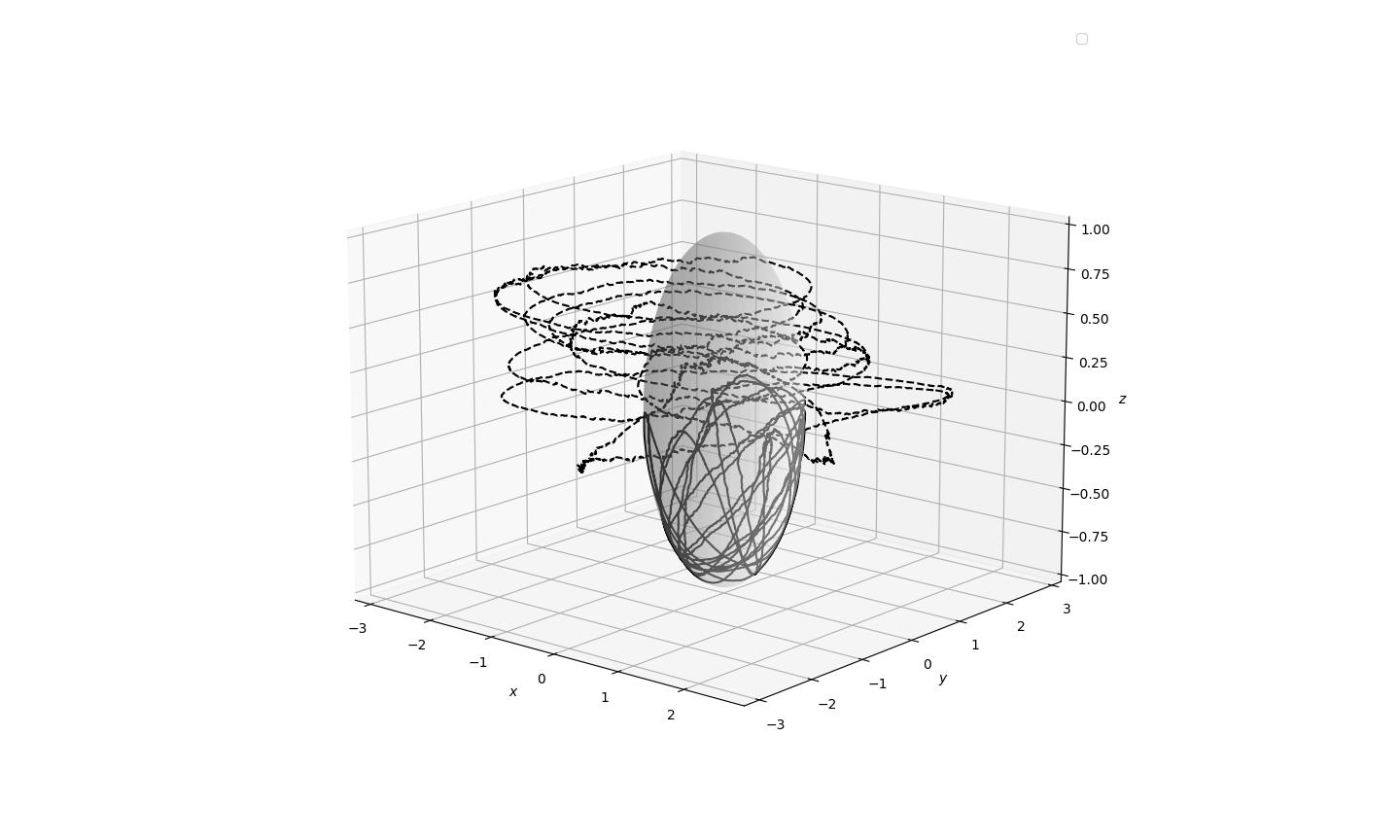}
\captionof{figure}{A single realisation of the stochastic heavy top generated with the TMK method with $\mathbb{I}={\rm diag}(4,4,1)$. The angular momentum (dashed) $\boldsymbol \pi$ wanders around the sphere and the gravity vector $\boldsymbol \gamma$ exactly stays on the sphere.}
\label{fig:dynamics2}
\end{minipage}
\begin{minipage}{.02\textwidth}
\hfill
\end{minipage}
\begin{minipage}[t]{.47\textwidth}
\vspace{-2ex}
\includegraphics[width=.99\textwidth]{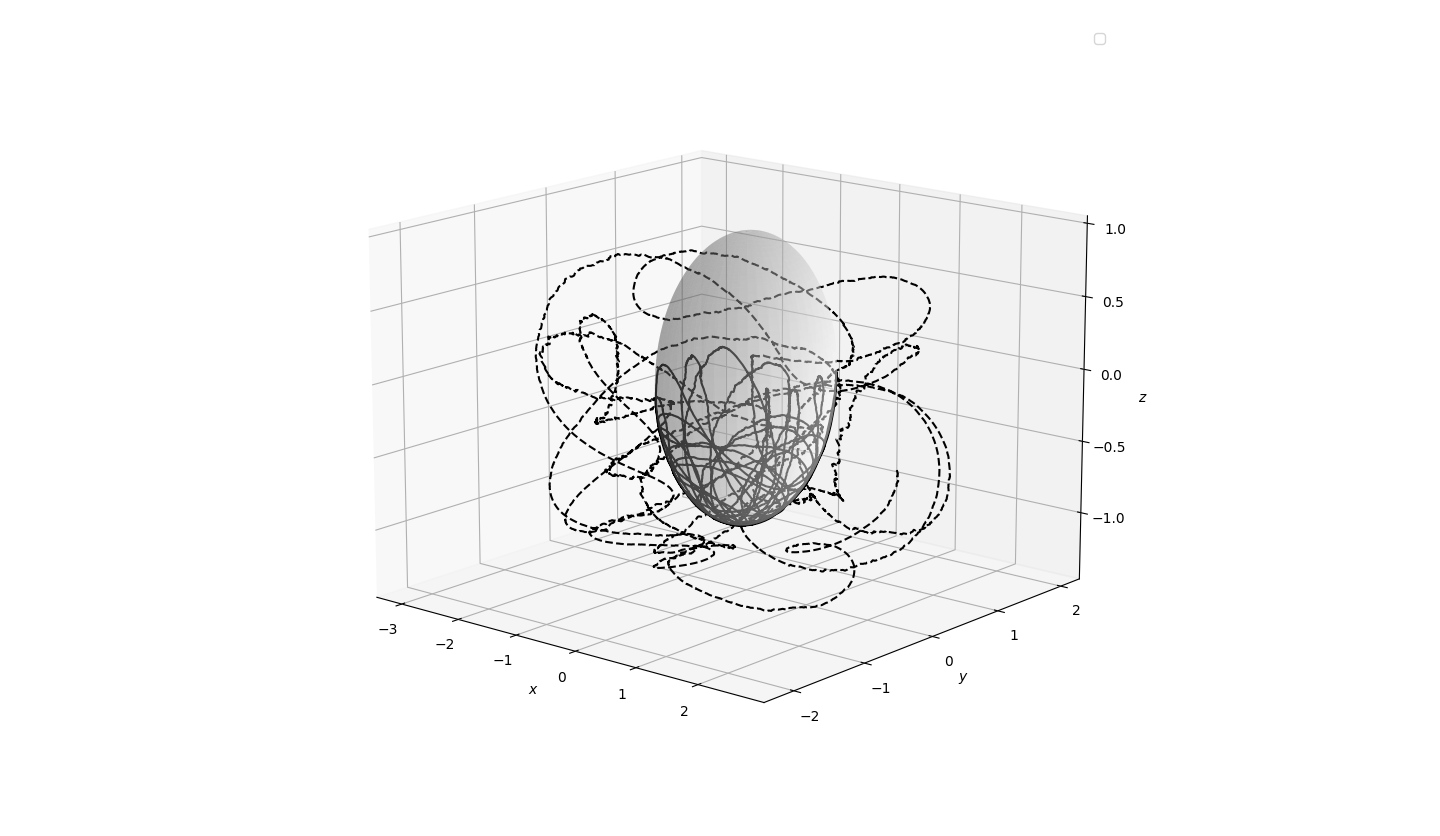}
\captionof{figure}{A single realisation of the stochastic heavy top generated with the TMK method with $\mathbb{I}={\rm diag}(4,2,1)$. The angular momentum (dashed) $\boldsymbol \pi$ wanders around the sphere and the gravity vector $\boldsymbol \gamma$ exactly stays on the sphere.}
\label{fig:dynamics1}
\end{minipage}
\medskip

In figures \ref{fig:casimir1} and \ref{fig:casimir2} it is evident that the IM is not able to conserve the Casimirs exactly, in contrast to the TMK method. In figure \ref{fig:casimirsIMMK} the relative error of both Casimirs is plotted for the TMK method, which shows that TMK is able to represent conservation of Casimirs to machine accuracy. The linear trend that is visible in the evolution of the Casimirs for the TMK method in figure \ref{fig:casimirsIMMK} is a result of accumulation of round-off errors. Extrapolating this linear growth would mean that after approximately $10^{10}$ time steps, the error in the TMK method would be comparable to the error of the IM method. Upon removing the noise, by setting $\boldsymbol \alpha = \mathbf{0}$, the TMK method gains conservation of the Hamiltonian as an additional property. This is shown in figure \ref{fig:energy}.
\medskip

\noindent
\begin{minipage}[t]{.47\textwidth}
\vspace{-2ex}
\includegraphics[width=.99\textwidth]{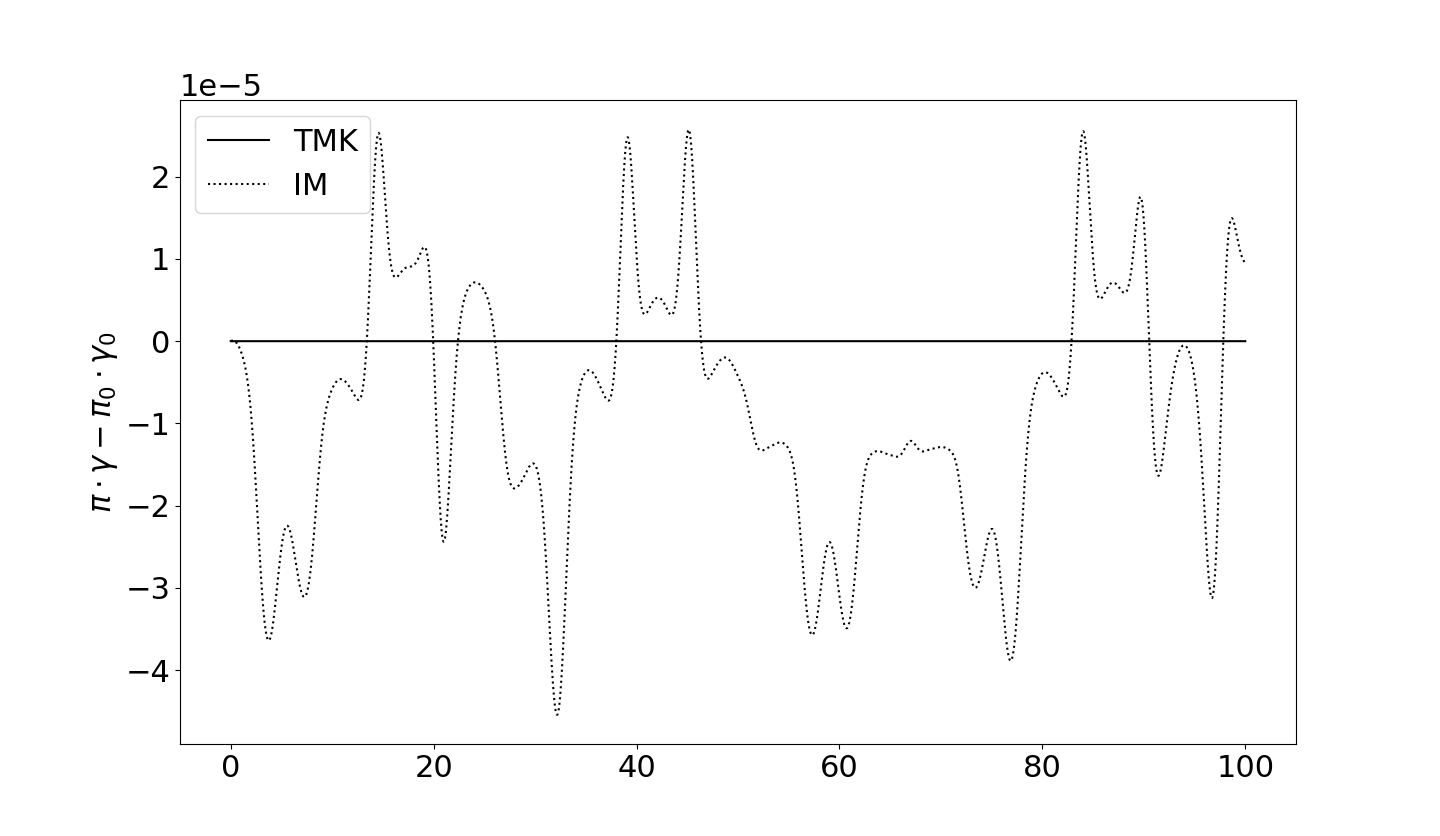}
\captionof{figure}{A single realisation of the Casimir $\boldsymbol \pi\cdot\boldsymbol \gamma$ compared to its initial value. Results generated by the IM method applied directly to the Lie-Poisson equation and by the TMK method are compared.}
\label{fig:casimir1}
\end{minipage}
\begin{minipage}{.02\textwidth}
\hfill
\end{minipage}
\begin{minipage}[t]{.47\textwidth}
\vspace{-2ex}
\includegraphics[width=.99\textwidth]{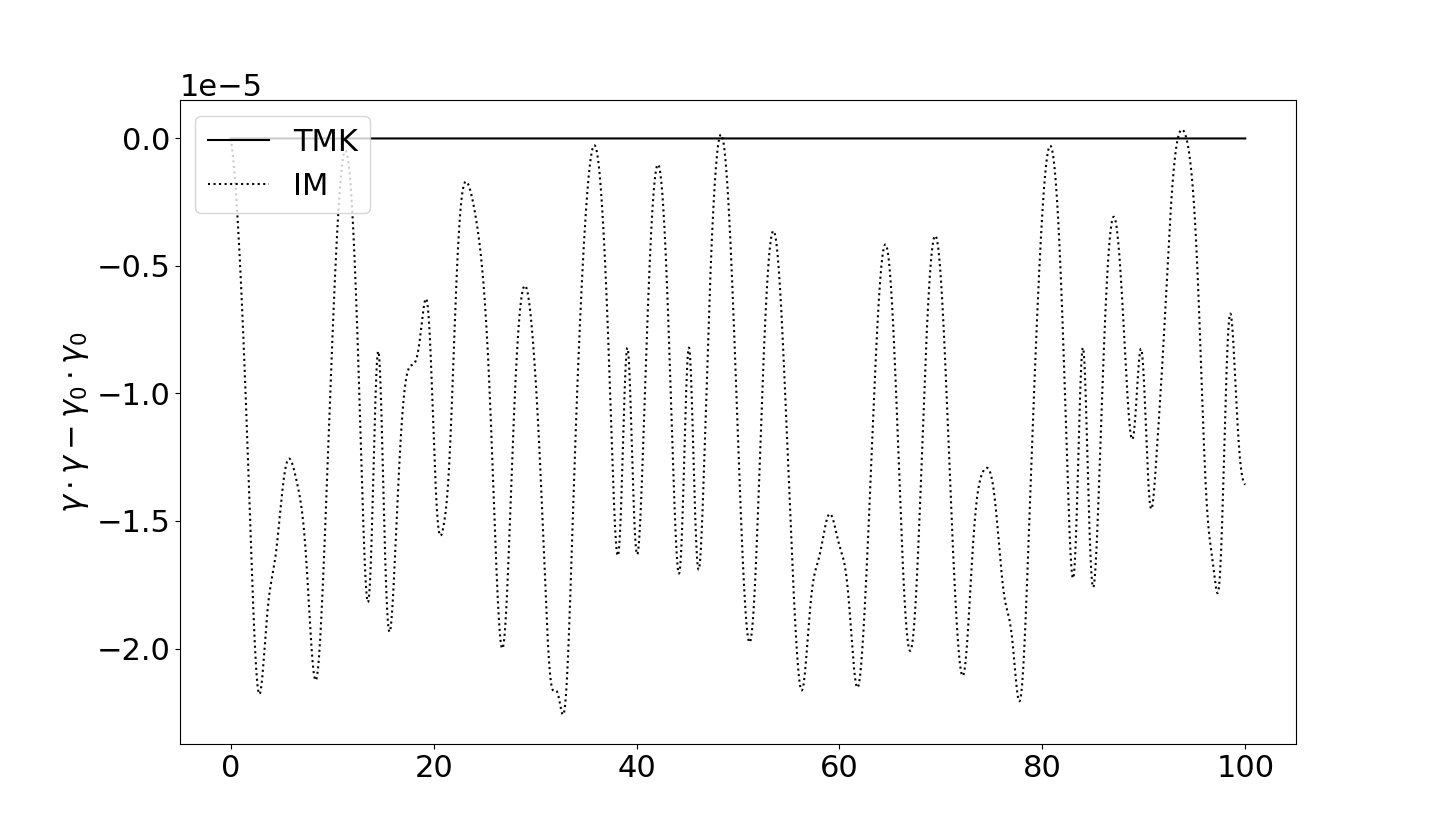}
\captionof{figure}{A single realisation of the Casimir $|\boldsymbol \gamma|^2$ compared to the initial value. Results generated by the IM method applied directly to the Lie-Poisson equation and by the TMK method are compared.}
\label{fig:casimir2}
\end{minipage}
\medskip

\noindent
\begin{minipage}[t]{.47\textwidth}
\vspace{-2ex}
\includegraphics[width=.99\textwidth]{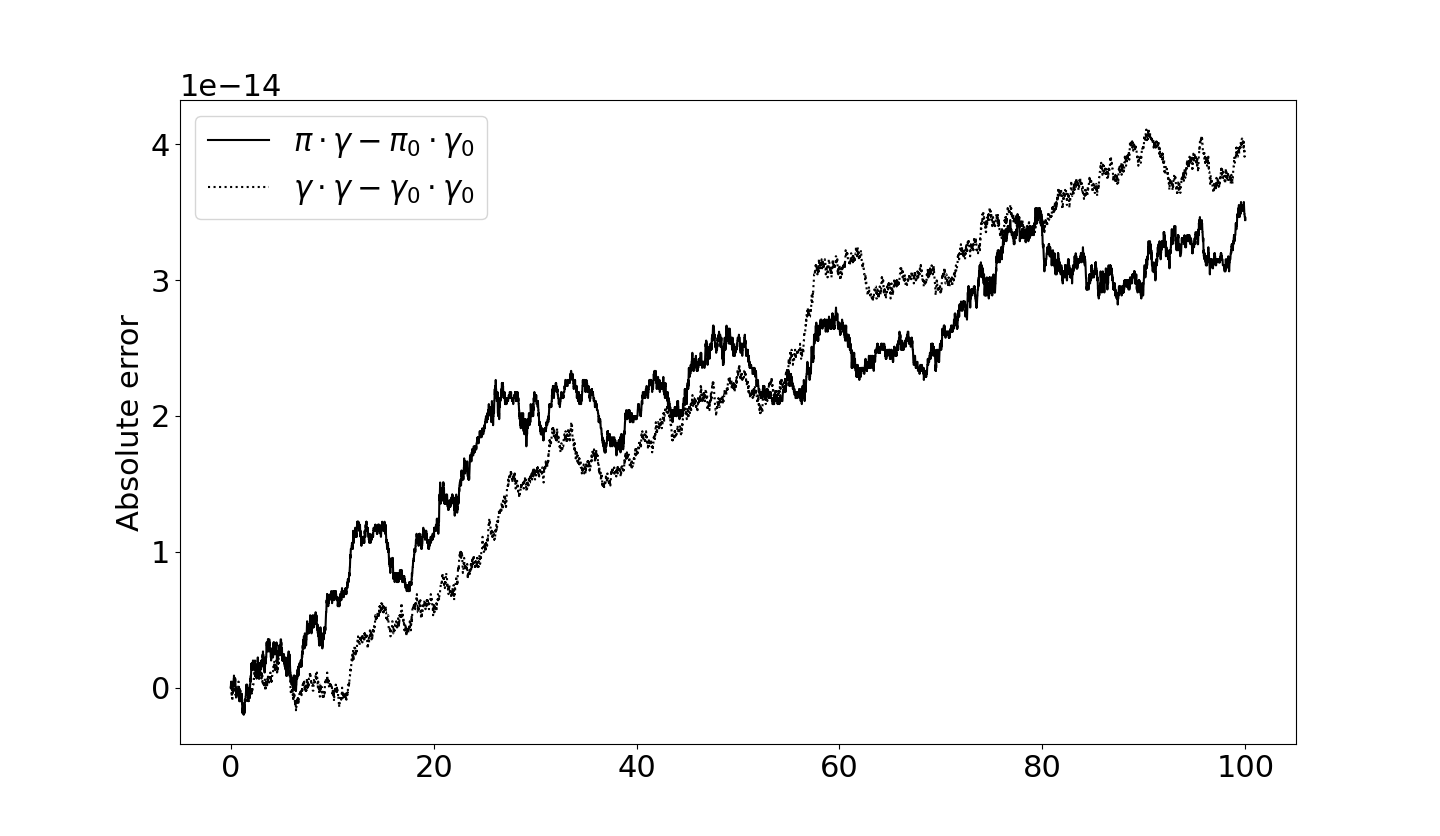}
\captionof{figure}{A single realisation of the absolute error of both Casimirs generated by the TMK method.}
\label{fig:casimirsIMMK}
\end{minipage}
\begin{minipage}{.02\textwidth}
\hfill
\end{minipage}
\begin{minipage}[t]{.47\textwidth}
\vspace{-2ex}
\includegraphics[width=.99\textwidth]{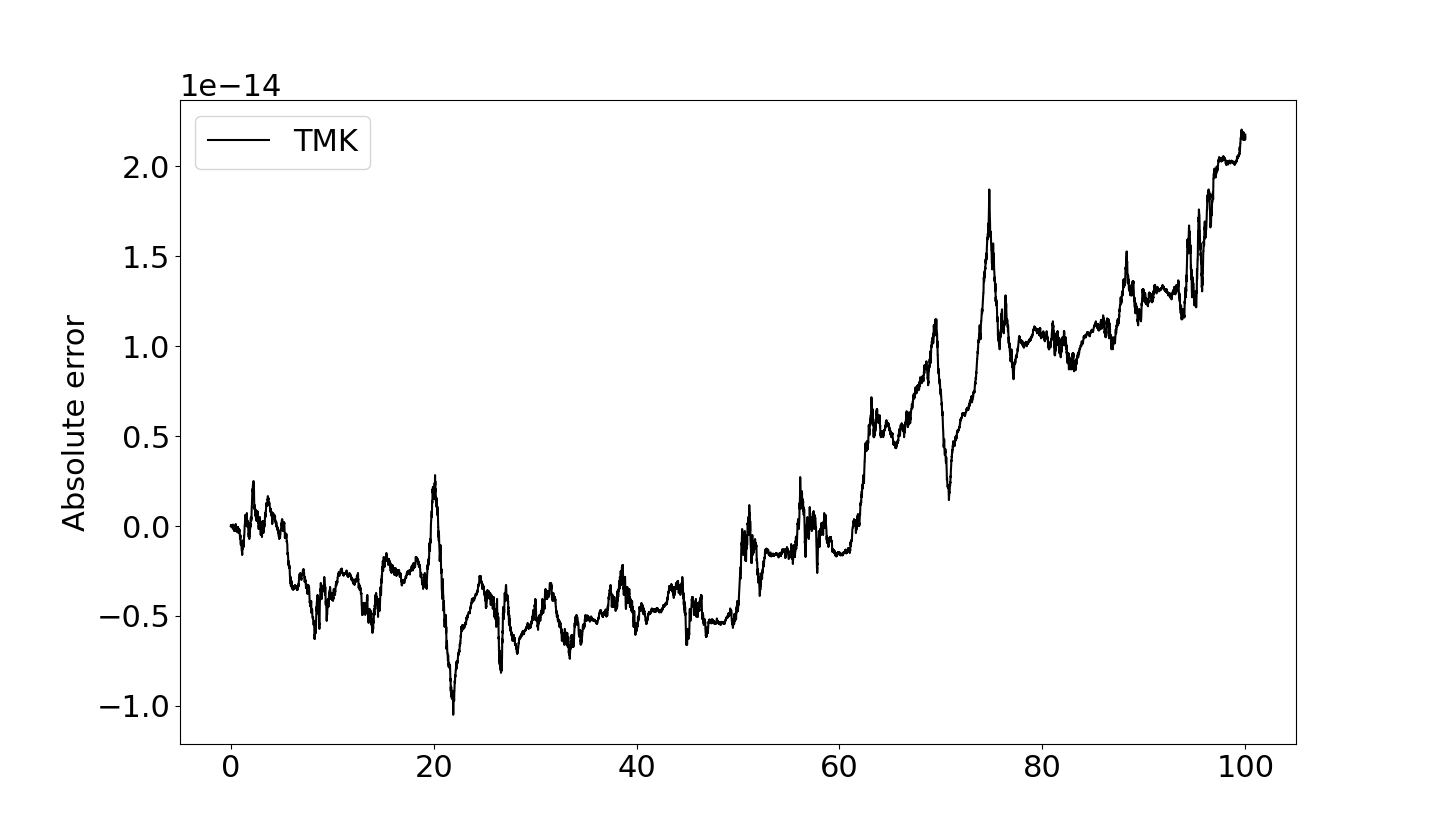}
\captionof{figure}{Upon setting $\boldsymbol \alpha=\mathbf{0}$, the noise is removed. In absence of noise, the Hamiltonian $\hslash$ of the heavy top is conserved by the TMK method. }
\label{fig:energy}
\end{minipage}
\medskip

\section{Sine-Euler}\label{sec:sEuler}
In this section the stochastic LP integrator based on the TMK method is illustrated for the sine-Euler equations introduced in \cite{zeitlin1991finite}. For completeness, we give a brief introduction of this fluid-mechanical model, after which the stochastic extension is studied for low-dimensional discrete truncation.

Consider the motion of an ideal incompressible fluid on the flat torus $\mathbb{T}^2$ governed by Euler's equations, which written for the vorticity $\omega$ read
\begin{equation}
\dot{\omega} = \left\{ \frac{\delta \hslash}{\delta \omega},\omega \right\},
\label{eq:Euler_vort}
\end{equation}
where $\{\cdot,\cdot\}$ is the Poisson bracket, defined as
\begin{equation}
\begin{aligned}
\{ f,g \} = \partial_x f \partial_y g - \partial_y f \partial_x g,  \hspace{3mm}\forall f,g \in C^{\infty} (\mathbb{T}^2),
\end{aligned}
\end{equation}
and $\hslash$ is the Hamiltonian
\begin{equation}
\hslash(\omega) = -\frac{1}{2} \int_{\mathbb{T}^2} \omega \psi d\Omega,
\label{eq:H_seuler}
\end{equation}
with $\psi$ the stream function, linked to the vorticity through $\Delta\psi=\omega$. System (\ref{eq:Euler_vort}) is a Lie-Poisson system on $C^{\infty} (\mathbb{T}^2)$ which admits infinitely many Casimir functions
\begin{equation}
\begin{aligned}
\mathcal{C}_k(\omega) = \int_{\mathbb{T}^2} \omega^k d\Omega, && k=1,2,\ldots
\end{aligned}
\label{eq:Casimirs_seuler}
\end{equation}
i.e, the integrated powers of vorticity are invariants of motion. In Fourier space equations (\ref{eq:Euler_vort}) become
\begin{equation}
\dot{\omega}_m = \sum_{n \neq 0} \frac{m \wedge n}{|n|^2} \omega_{m+n}\omega_{-n},
\label{eq:Euler_vort_fs}
\end{equation}
where $m=(m_1,m_2)$ is an integer vector and $m \wedge n = m_1n_2 - m_2n_1$. Numerical simulation of (\ref{eq:Euler_vort_fs}) requires a truncation to a finite set of Fourier modes. The strict conservation of the Casimirs is then no longer maintained. In fact, in order to preserve the underlying geometric structure also upon truncation at arbitrary finite order, \cite{zeitlin1991finite} put forward an approach based on the theory of geometric quantisation of \cite{hoppe1989diffeomorphism}. In this context, quantisation refers to the process of constructing a Lie algebra of $N \times N$ complex matrices, which replaces the Poisson bracket by the matrix commutator. In particular, there exists a basis of $\mathfrak{su}(N)$ (skew-Hermitian traceless matrices) with structure constants converging to those of the Fourier basis of $C^{\infty} (\mathbb{T}^2)$. In this basis, one can rewrite (\ref{eq:Euler_vort}) in terms of the vorticity matrix $W$:
\begin{equation}
\dot{W} = [ P,W ],
\label{eq:Euler_vort_mat}
\end{equation}
where $W, P \in \mathfrak{su}(N)$ with $P$ the stream matrix. As pointed out in \cite{zeitlin1991finite}, traces of powers of $W$,
\begin{equation}
\begin{aligned}
C_k(W)=\text{Tr}(W^k) && \text{for } k=1,\ldots,N,
\end{aligned}
\label{eq:power_vort}
\end{equation}
are conserved by (\ref{eq:Euler_vort_mat}). This is the discrete analogue of conservation of integrated powers of vorticity in the continuum. In Fourier coordinates, (\ref{eq:Euler_vort_mat}) gives the sine-Euler equations:
\begin{equation}
\dot{\omega}_m = \sum_{n = -K}^K \frac{1}{\varepsilon} \frac{\sin( \varepsilon m \wedge n)}{|n|^2} \omega_{m+n}\omega_{-n},
\label{eq:sEuler_vort_fs}
\end{equation}
with $K=(N-1)/2$ and $\varepsilon=2\pi/N$.

Following the SALT approach, we derive the stochastic Euler equations by introducing diffusion Hamiltonians
\begin{equation}
\hslash_s \circ d\mathbf{S}_t = -\frac{1}{2} \int_{\mathbb{T}^2} \psi \omega d\Omega dt - \frac{1}{2} \sum_{i=1}^M \int_{\mathbb{T}^2} \omega \zeta_i d\Omega \circ dW_t^i,
\label{eq:seuler_H_stoc}
\end{equation}
where $\zeta_i \in C^{\infty} (\mathbb{T}^2)$. Using $\hslash_s \circ d\mathbf{S}_t$ in (\ref{eq:Euler_vort}) and applying the sine function as in (\ref{eq:sEuler_vort_fs}) one arrives at the stochastic sine-Euler equations in Fourier space
\begin{equation}
d\omega_m = \sum_{n = -K}^K \frac{1}{\varepsilon} \sin( \varepsilon m \wedge n) \omega_{m+n} \left( \frac{\omega_{-n}}{|n|^2} dt + \sum_{i=1}^M \zeta_{i,-n} \circ dW_t^i \right).
\label{eq:sEuler_vort_fs_stoc}
\end{equation}

To illustrate the TMK method specified in sections \ref{sec:slp} and \ref{sec:slpdynamics} for the sine-Euler equations, we simulate the dynamics for $N=3$. Differently from the heavy top (section \ref{sec:ht}), this system has a quadratic as well as a cubic Casimir, being $\mathcal{C}_2=\text{Tr}(W^2)$ and $\mathcal{C}_3=\text{Tr}(W^3)$, respectively. By simulating the sine-Euler equations for $N>2$ we demonstrate conservation of polynomial invariants of high order.

As initial condition, we set the Fourier coefficients of vorticity to be randomly distributed over a uniform distribution between $-0.5$ and $0.5$. The time step is $\Delta t=0.5$ and the final simulation time is $T=10000$ to reach a well-developed statistically stationary solution. Noise is injected at modes $n=\left[(1,1),(1,-1)\right]$ with amplitude $10^{-1}$, representing a high-frequency disturbance. Each mode has its own independent Wiener process. At any time $t$, the state of the system is given by the four-dimensional state vector $(\omega_{0,1},\omega_{1,1},\omega_{1,0},\omega_{1,-1})$. Conservation of the Casimirs is shown in figure \ref{fig:seuler_casimirs} and the Casimir behaviour corresponding to direct integration of the sine-Euler equations using the trapezoidal rule is shown in figure \ref{fig:seuler_casimirs_T}. For the sine-Euler equations, the TMK method was empirically found to be stable even for time step sizes beyond $\Delta t =500$. This is far beyond the time step size where the direct trapezoidal rule breaks down, which was empirically found to be around $\Delta t=2$. For a fair comparison, we choose $\Delta t=0.5$, since both methods are stable.
\medskip

\noindent
\begin{minipage}[t]{.47\textwidth}
\vspace{-2ex}
\includegraphics[width=.99\textwidth]{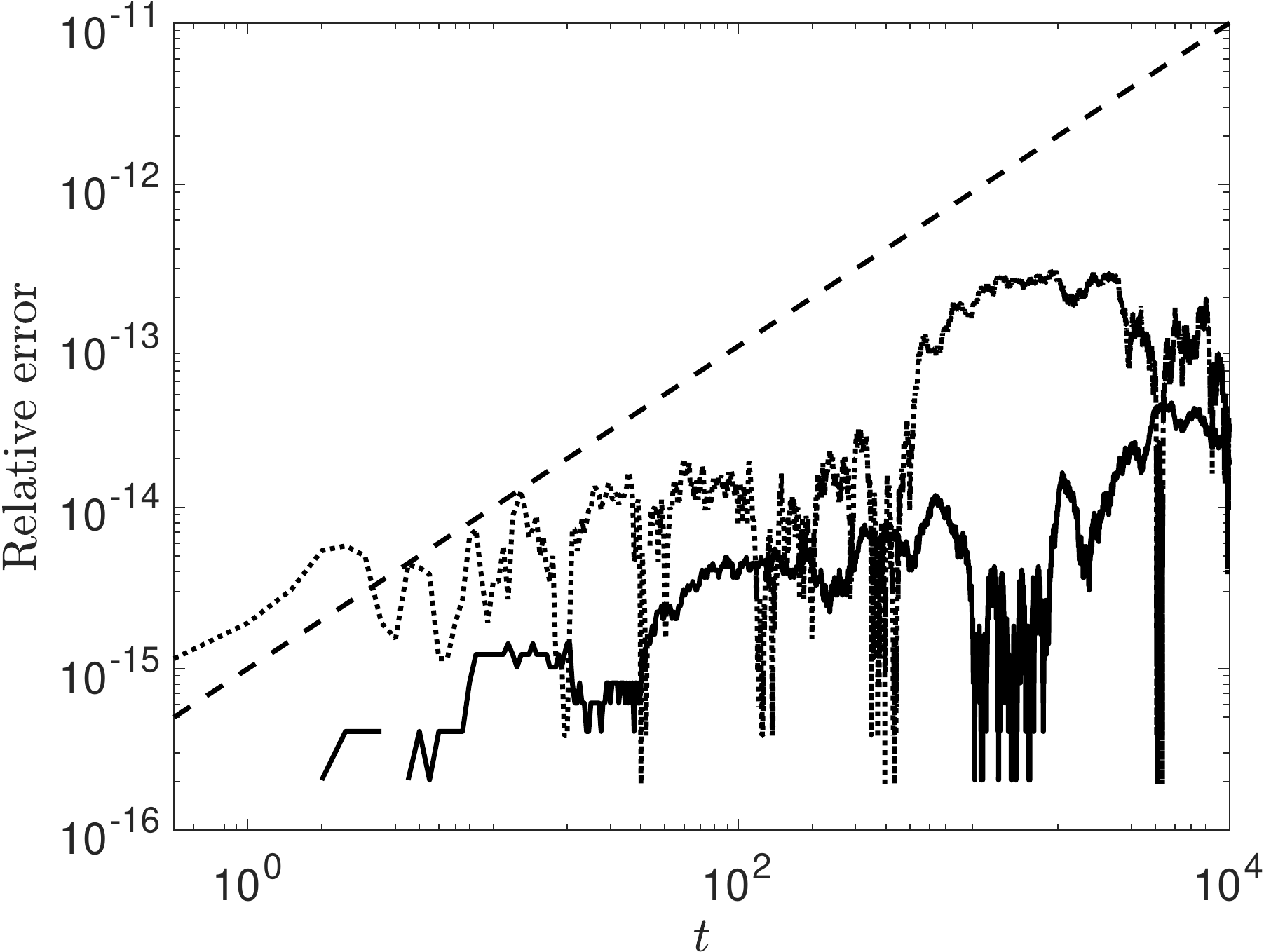}
\captionof{figure}{Relative error of the Casimirs for one stochastic realization of the sine-Euler equations simulated by using the TMK method. The solid line refers to the conservation of $C_2$, while the dash-dotted line refers to the conservation of $C_3$.}
\label{fig:seuler_casimirs}
\end{minipage}
\begin{minipage}{.02\textwidth}
\hfill
\end{minipage}
\begin{minipage}[t]{.47\textwidth}
\vspace{-2ex}
\includegraphics[width=.99\textwidth]{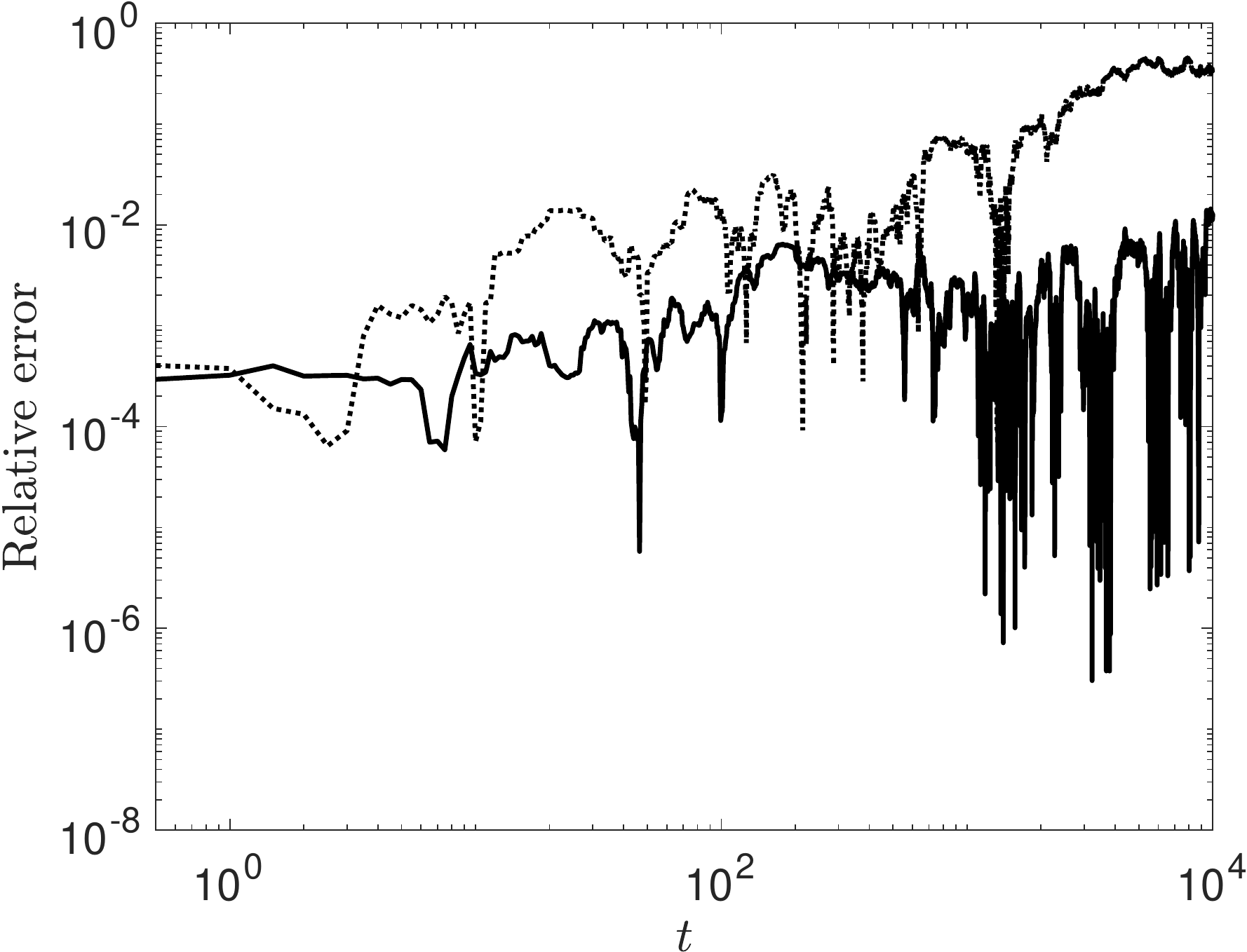}
\captionof{figure}{Relative error of the Casimirs for one stochastic realization of the sine-Euler equations simulated by using the trapezoidal rule applied directly to the LP equations. The solid line refers to the conservation of $C_2$, while the dash-dotted line refers to the conservation of $C_3$.}
\label{fig:seuler_casimirs_T}
\end{minipage}

By construction the TMK integrator conserves the Casimirs up to machine precision. Furthermore, analogously to the heavy top, the relative error appears to be robust and increase approximately linearly, as indicated by the dashed line in figure \ref{fig:seuler_casimirs}, over a rather long simulation time.

\section{Conclusion}\label{sec:conclusion}
In this paper we introduced numerical methods for stochastic Lie-Poisson (LP) equations that preserve the coadjoint orbit structure by extending the Runge-Kutta-Munthe-Kaas (RKMK) method to stochastic differential equations (SDEs). Specifically, this method is able to preserve Casimirs exactly. The stochasticity is defined with respect to the Stratonovich integral, since the ordinary chain rule is required to rigorously extend the approach to stochastic dynamics. The noise was chosen to be of the stochastic advection by Lie transport (SALT) type, i.e., multiplicative noise coupled linearly to the momentum variable. Moreover, given that the Stratonovich integral is defined at the midpoint of the integrand, the implicit midpoint rule is a natural choice for the integration of stochastic LP equations with Stratonovich noise. We applied the Munthe-Kaas approach to obtain a stochastic differential equation on the Lie algebra which we solve using the implicit midpoint rule. The implicit midpoint method is in the class of RKMK methods. The solution of the SDE on the Lie algebra is used to generate a group element that together with the coadjoint representation of the group on the dual of the Lie algebra generates the coadjoint orbit associated with the initial condition. Exactly generating coadjoint orbits guarantees the conservation of Casimirs. The implicit midpoint rule is particularly convenient, because in absence of the noise, the integrator preserves the deterministic Hamiltonian. 

The implied SALT-induced SDE on the Lie algebra has additive noise, whereas the stochastic LP equation on the dual of the Lie algebra has multiplicative noise. To illustrate the theoretical results with numerical experiments, we applied the implicit midpoint (IM) rule to the LP equations directly and used the trapezoidal rule (TMK) within the class of stochastic LP integrators. A comparison of IM and TMK showed that TMK performs according to theory when conservation of Casimirs is concerned, apart from small trends in the error due to round-off effects. This was considered for the examples of the heavy top and for low-order truncation of the sine-Euler equations. The numerical illustrations clarify that the implementation of the numerical integrators is fully inline with the theoretical preservation properties. By switching off the noise and using TMK to solve for the deterministic heavy top, we also showed that TMK conserves the Hamiltonian.

The developed stochastic LP integrator is invaluable for the long-time simulation of stochastic mechanical systems for which the conservation of the geometric structure is essential. Apart from the two test-cases, i.e., the heavy top and the sine-Euler equations, the new LP integrator can be applied effectively to a range of applications, among which are robotics, celestial mechanics, biomechanics, rigid body mechancis etc. These LP mechanical systems are perturbed stochastically with the SALT method from \cite{holm2015variational}, which also enables a basis for uncertainty quantification through stochastic forcing. An extension to stochastic affine LP equations, which arise in machine learning and mechanics on centrally extended Lie algebras is left for future work.

\section*{Abbreviations}
\begin{table}[H]
\centering
\begin{tabular}{l|l}
Abbreviation & Meaning\\
\hline
LP & Lie-Poisson\\
RKMK & Runge-Kutta-Munthe--Kaas\\
IM & Implicit midpoint\\
TMK & Trapezoidal Munthe-Kaas\\
SDE & Stochastic differential equation\\
SALT & Stochastic advection by Lie transport
\end{tabular}
\end{table}

\section*{Availability of data and material}
No data was used in the preparation of this work.

\section*{Competing interests}
The authors have no competing interests to declare that are relevant to the content of this work. All authors certify that they have no affiliations with or involvement in any organisation or entity with any financial interest or non-financial interest in the subject matter or materials discussed in this manuscript.

\section*{Funding}
This work was performed in the project SPReStO (structure-preserving regularization and stochastic forcing for nonlinear hyperbolic partial differential equations), supported by a NWO TOP 1 grant.

\section*{Author's contributions}
All authors contributed equally to this work.

\section*{Acknowledgements}
We wish to express our gratitude towards Raffaele D'Ambrosio, Arnout Franken, Darryl Holm and Ruiao Hu for their valuable input in the preparation of this paper. We also wish to acknowledge the anonymous reviewer for their valuable comments. EL and SE wish to thank the Gran Sasso Science Institute for its warm hospitality during a visit in 2021. This work was performed in the project SPReStO (structure-preserving regularization and stochastic forcing for nonlinear hyperbolic PDEs), supported by a NWO TOP 1 grant.

\bibliographystyle{plainnat}      
\bibliography{biblio}

\end{document}